\RequirePackage{fix-cm}
\pdfoutput=1

\documentclass[leqno, fleqn, centertags, 12pt]{article}

\usepackage{latexsym}
\usepackage{amsmath}
\usepackage{amssymb}
\usepackage{verbatim}

\usepackage[T1]{fontenc}

\usepackage[upright, widespace]{fourier}

\usepackage{fullpage}

\usepackage{tikz-cd}

\tikzcdset{row sep/hugeplus/.initial=6em}
\tikzcdset{row sep/hugepl/.initial=5.4em}
\tikzcdset{row sep/special/.initial=12ex}
\tikzcdset{row sep/spe/.initial=10.5ex}
\tikzcdset{column sep/spec/.initial=5em}
\tikzcdset{column sep/specc/.initial=8em}
\tikzcdset{column sep/speccc/.initial=10em}

\usepackage{comment}

\newskip\nineskipamount \nineskipamount=9pt plus 0pt minus 0pt
\newskip\zeroskipamount \zeroskipamount=0pt plus 0pt minus 0pt
\usepackage[noorphans,vskip=\nineskipamount]{quoting}

\makeatletter
\renewcommand{\@makefntext}[1]{\vspace*{0.5ex}\parindent=0em
\hspace*{-0.4em}
\hbox to 0.4em{\hss\@makefnmark}\hspace*{0.4em}{#1}
}
\makeatother

\newcounter{mysectionnumber}
\newcommand{\mysection}[2]{\setcounter{footnote}{0}
\setcounter{myparnum}{0}
\refstepcounter{mysectionnumber}
\vspace{27pt}{\Large {\themysectionnumber.} {#1}}\label{#2}\vspace*{15pt}}

\newcommand{\myit}[1]{\textbf{\textit{#1}}\hspace{0.0em}}

\newcounter{myparnum}[mysectionnumber]
\renewcommand{\themyparnum}{\arabic{mysectionnumber}.\arabic{myparnum}}
\newcommand{\mypar}[2]{\refstepcounter{myparnum}{\vspace{\medskipamount}\textbf{{\themyparnum. #1}\label{#2}}\hspace{0.5em}}}

\newcounter{mylemmanum}[myparnum]
\newcommand{\myuppar}[1]{\vspace{\medskipamount}\textbf{#1}\hspace*{0.5em}}

\newcounter{myappendnumber}
\newcounter{myaparnum}[myappendnumber]
\newcounter{myapparnum}[mysectionnumber]
\newcommand{\proof}{\vspace{\medskipamount}{\textbf{{\emph{Proof}.}}\hspace*{1em}}}

\newcommand{\eproof}{ $\blacksquare$}

\newcommand{\dis}{\displaystyle}

\def\sss{\hspace{0.05em}\ }
\def\dss{\hspace{0.1em}\ }
\def\trs{\hspace{0.15em}\ }
\def\qss{\hspace{0.2em}\ }
\def\pss{\hspace{0.3em}\ }
\def\oss{\hspace{0.4em}\ }

\def\halfff{\hspace*{0.025em}}
\def\fff{\hspace*{0.05em}}
\def\dff{\hspace*{0.1em}}
\def\trf{\hspace*{0.15em}}
\def\qff{\hspace*{0.2em}}
\def\pff{\hspace*{0.3em}}
\def\off{\hspace*{0.4em}}

\newcommand{\nsp}{\hspace*{-0.1em}}
\newcommand{\nnsp}{\hspace*{-0.15em}}

\newcommand{\dnsp}{\hspace*{-0.2em}}

\renewcommand{\leq}{\leqslant}
\renewcommand{\geq}{\geqslant}

\newcommand{\zzz}{\mathbb{Z}}

\newcommand{\hclass}[1]{[\dff #1 \dff]}

\newcommand{\ttoo}{\hspace*{0.2em}\longrightarrow\hspace*{0.2em}}

\renewcommand{\hom}{\operatorname{Hom}\trf}

\begin{document}

\setlength{\baselineskip}{12pt plus 0pt minus 0pt}
\setlength{\parskip}{12pt plus 0pt minus 0pt}
\setlength{\abovedisplayskip}{12pt plus 0pt minus 0pt}
\setlength{\belowdisplayskip}{12pt plus 0pt minus 0pt}

\newskip\smallskipamount \smallskipamount=3pt plus 0pt minus 0pt
\newskip\medskipamount   \medskipamount  =6pt plus 0pt minus 0pt
\newskip\bigskipamount   \bigskipamount =12pt plus 0pt minus 0pt

\author{Nikolai\qss V.\qss Ivanov}
\title{Non-abelian\qss cohomology\qss and\qss Seifert--van~Kampen\qss theorem}
\date{}

\footnotetext{\hspace*{-0.65em}\copyright\oss 
Nikolai\qss V.\qss Ivanov,\oss 2023.\trs}

\maketitle

{
\renewcommand{\baselinestretch}{1}
\selectfont

\myit{\hspace*{0em}\large Contents}  \vspace*{1.5ex}\\ 
\hbox to 0.8\textwidth{\myit{1.}\hspace*{0.5em} Introduction \hfil 1}  \hspace*{0.5em} \vspace*{0.25ex}\\
\hbox to 0.8\textwidth{\myit{2.}\hspace*{0.5em} Group-valued cochains in dimensions $0$ and $1$ \hfil 6}  \hspace*{0.5em} \vspace*{0.25ex}\\
\hbox to 0.8\textwidth{\myit{3.}\hspace*{0.5em} Short\sss cochains \hfil 8}  \hspace*{0.5em} \vspace*{0.25ex}\\
\hbox to 0.8\textwidth{\myit{4.}\hspace*{0.5em} The standard\dss Seifert--van Kampen\dss theorem \hfil 11}  \hspace*{0.5em} \vspace*{0.25ex}\\
\hbox to 0.8\textwidth{\myit{5.}\hspace*{0.5em} Unions of\dss several\sss subsets \hfil 15}
\vspace*{0.25ex}\\
\hbox to 0.8\textwidth{\myit{6.}\hspace*{0.5em} van Kampen\dss theorems \hfil 17}  \hspace*{0.5em} \vspace*{1.5ex}\\
\hbox to 0.8\textwidth{\myit{References}\hspace*{0.5em} \hfil 26} \hspace*{0.5em}  \vspace*{0.25ex}

\renewcommand{\baselinestretch}{1}
\selectfont

\mysection{Introduction}{introduction}\vspace{-0.25pt}

\myuppar{The standard\dss Seifert--van Kampen\dss theorem.}
The\dss Seifert--van Kampen\dss theorem\dss is\dss a common name for\sss
theorems relating\sss the fundamental\sss group 
of\dss the union of\dss two or more spaces\sss to\sss the fundamental groups of these spaces
and\sss their pairwise intersections,\qss
under suitable\sss local\sss conditions.\oss
The most\sss familiar\dss Seifert--van Kampen\dss theorem\dss
deals with a space $X$ presented as\sss the union of\dss two open subsets $U\fff,\qff V$
such\sss that\sss $X\fff,\qff U\fff,\qff V$ and\sss the intersection\sss $U\dff \cap\dff V$
are path-connected and asserts\sss that\sss
$\pi_{\dff 1}\dff(\trf X\trf)$\sss is\dss equal\dss to\sss the free product\sss
of\dss $\pi_{\dff 1}\dff(\trf U\trf)$\sss and\sss $\pi_{\dff 1}\dff(\trf V\trf)$\sss
amalgamated over\sss $\pi_{\dff 1}\dff(\trf U\dff \cap\dff V\trf)$\nnsp.\oss
The assumption\sss that\sss the subsets $U\fff,\qff V$ are open was,\oss
probably,\oss first\sss introduced\dss by\dss R.\dss Crowell\sss and\dss R.\dss Fox\qss 
\cite{c},\trs \cite{cf}.\oss
Seifert\dss and\dss van Kampen\dss worked\sss with closed subsets,\oss
and\dss Seifert\dss even with subcomplexes of\dss triangulated spaces.\oss
In applications\sss the subsets are usually closed,\qss but\sss
are deformation retracts of\dss their open neighborhoods\sss
and can be replaced\dss by\sss such\sss neighborhoods.\oss
At\sss the same\sss time\sss dealing\sss first\sss with open subsets,\qss
or subsets with\sss the interiors covering $X$\nnsp,\qss
allows\sss to separate\sss the global\dss issues,\oss
present\sss already for open subsets,\oss
from\sss the\sss local\sss ones involved\sss in passing\sss from
open subsets\sss to closed ones.\qss 
The present\sss paper deals with\sss the global\dss issues 
and\sss mostly with open subsets.\vspace{-1pt}

R.\dss Fox\dss related\qss \cite{f}\qss that\sss he 
introduced\dss the name\dss ``van Kampen\dss theorem''\dss
for\sss the special\sss case of\dss 
path-connected\sss $U\fff,\qff V$ and\sss $U\dff \cap\dff V$\dnsp,\oss
overlooking\sss the fact\sss that\sss this case was
proved\dss by\dss Seifert\qss \cite{s}\qss two years before\dss
van Kampen's\dss work\qss \cite{vk}.\oss
Unfortunately,\oss this name\dss is\dss still\dss widely used\sss in\sss same sense.\oss
This\sss led\dss to many claims\sss
that\dss van Kampen\dss theorem\dss is\dss not\sss sufficient\sss
to compute even\sss the fundamental\dss group of\dss the circle,\qss
because\sss the circle cannot\sss be presented as\sss the union
of\dss two path-connected subsets with path-connected\sss intersection.\oss
In\sss fact,\oss the computation of\dss the fundamental\sss group of\dss 
the circle\dss is\dss a very special\sss case of\trs van Kampen\dss results.\vspace{-1pt}

\myuppar{van Kampen\dss theorems.}
Let\sss $C\fff,\pff B$ be\sss topological\sss spaces
such\sss that\sss $B$\sss is\dss path-con\-nect\-ed and $C$ consists
of\dss finite or countable number of\dss path-components.\oss
Let $B_{\dff 1}\dff,\qff B_{\dff 2}\dff,\qff \ldots$ be a finite or countably\sss infinite
collection of\dss closed subsets of\sss $C$ homeomorphic\sss to $B$ and such\sss that\sss
every\sss path-connected component\sss of\sss $C$ contains at\sss least\sss
one set\sss $B_{\dff i}$\nsp.\oss
Suppose\sss that\sss some homeomorphisms\dss 
$h_{\dff i}\dff \colon\dff B_{\dff i}\qff \ttoo\qff B$\dss
are fixed,\oss
and\dss let\dss $A$\dss be\sss the quotient\sss space of\dss $C$\dss obtained\dss by\sss
identifying\dss every\dss $B_{\dff i}$\dss with\dss $B$\dss 
by\sss $h_{\dff i}$\nsp.\qss
van Kampen\qss \cite{vk}\qss described\sss $\pi_{\dff 1}\dff(\trf A\trf)$\sss 
in\dss terms of\sss $\pi_{\dff 1}\dff(\trf C\trf)$\sss and\sss $\pi_{\dff 1}\dff(\trf B\trf)$\sss
and\dss of\dss homeomorphisms\dss $h_{\dff i}$\dss under some\sss local\sss assumptions
about\sss subspaces\dss $B_{\dff i}$ 
and some\dss ``niceness''\dss assumptions about\sss the\sss topological\sss spaces involved.\oss
It\sss seems\sss that\sss the\sss latter are stronger\dss
than necessary\sss and were borrowed\dss from\dss Lefschetz's\dss
treatment\dss of\trs the excision\dss property\sss of\dss homology\sss groups.\oss
In\sss fact,\oss van Kampen\dss had no notion of\dss a quotient\sss space at\sss his
disposal,\oss and\sss used a somewhat\sss cumbersome description
of\dss relations between\sss $A$\nnsp,\dss $C$\nnsp,\oss
and\sss $B_{\dff i}$\nsp.\oss
He started with $A$ and\sss $B$ and constructed $C$\nnsp.\oss
In any case,\oss if\dss $C\off =\off [\dff 0\fff,\qff 1\dff]$\nnsp,\dss
$B_{\dff 1}\off =\off \{\trf 0\trf\}$\nnsp,\dss
$B_{\dff 2}\off =\off \{\trf 1\trf\}$\nnsp,\qss
than $A$ is\dss a circle and
$\pi_{\dff 1}\dff(\trf A\trf)\off =\off \zzz$ by\dss van Kampen\dss results.\oss\vspace{-1pt}

In\sss the\sss last\sss section of\dss van Kampen's\dss paper\qss
\emph{``the path\dss is\dss shown\sss to a more general\dss theorem,\oss
of\dss which however\sss the general\sss formulation\dss would\dss be more confusing\sss than helpful,\oss
so\sss that\dss it\dss is\dss suppressed''.}\oss
van Kampen\dss singles out\sss two most\sss important\sss special\sss cases of\dss his\sss theorem,\oss
both of\dss which are deal\sss with\sss the case of\dss 
two subsets $B_{\dff 1}\dff,\qff B_{\dff 2}$\nsp.\oss
In\sss the first\sss one\sss the space $C$\sss is\dss assumed\dss to be path-connected.\oss
In\sss the second special\sss case $C$\sss is\dss assumed\dss to
consist\sss of\dss two path-connected components 
$C_{\dff 1}\fff,\qff C_{\dff 2}$
containing $B_{\dff 1}\fff,\qff B_{\dff 2}$ respectively.\oss
See\qss \cite{vk},\oss Corollaries\qss 1\qss and\qss 2.\oss
The framework of\dss the second case\dss is\dss the same as\dss Seifert's\dss one.\oss
For a more detailed discussion of\trs van Kampen's\dss results we refer\sss to\dss
A.\dss Gramain\qss \cite{g}.\oss

van Kampen's\dss paper\qss \cite{vk}\qss has\sss the reputation of\dss being\sss difficult.\oss
It\dss is\dss indeed very densely written,\oss but\sss the present\sss author believes\sss
that\sss the reason\dss is\dss different.\oss
van Kampen\dss deals simultaneously with\sss three different\sss problems\fff:\oss
with\sss the\sss lack of\dss the notion of\dss quotient\sss spaces,\oss
with\sss the\sss local\sss issues caused\sss by working with closed subsets\qss
(for\sss the sake of\dss applications),\oss
and,\oss finally,\oss with algebraic issues associated nowadays with\sss the\sss term\dss
``van Kampen\dss theorem''.\oss

\myuppar{van Kampen's\qss framework and unions of\dss open subspaces.}
In contrast\sss with\dss Seifert\dss and\sss the modern expositions,\oss
van Kampen\dss worked not\sss with unions,\oss but\sss with some quotient\sss spaces.\oss
There\dss is\dss a simple\sss trick allowing simultaneously\sss to pass from\sss
van Kampen's\dss framework\sss to unions of\dss open subspaces
and\sss to separate\sss the\sss local\sss issues from\sss the global\sss ones.\oss
Let\dss\vspace{0.625pt}
\[
\quad 
B_{\dff \bullet}
\off =\off 
\bigcup\nolimits_{\trf i}\dff B_{\dff i}
\]

\vspace{-12pt}\vspace{0.625pt}
and\dss let\dss $h_{\dff \bullet}\dff \colon\dff B_{\dff \bullet}\qff \ttoo\qff B$\dss
be\sss the map defined\dss by\dss the maps\dss $h_{\dff i}$\nsp.\oss
Let\sss $X$\sss be\sss the cylinder of\trs the map\dss $h_{\dff \bullet}$\nsp.\oss
In more details,\qss $X$\sss is\dss the result\sss of\dss 
glueing of\dss the subset\sss\vspace{3pt}
\[
\quad 
C\qff \cup\qff B_{\dff \bullet}\dff \times\dff [\trf 0\fff,\qff 1\trf]
\qff \subset\qff C\dff \times\dff [\trf 0\fff,\qff 1\trf]
\]

\vspace{-9pt}
to $B$\sss by\sss the map\sss 
$B_{\dff \bullet}\dff \times\dff 1\dff \ttoo\qff B$\sss
induced\sss by $h_{\dff \bullet}$\nsp.\oss
There\dss is\dss an obvious map\dss $X\qff \ttoo\qff A$\nnsp.\oss
If\dss the pair\sss $(\trf C\fff,\pff B_{\dff \bullet}\trf)$\sss
has\sss the homotopy\sss extension property\halfff,\oss
i.e.\qss the inclusion\dss $B_{\dff \bullet}\qff \ttoo\qff C$\dss is\dss a cofibration,\oss
then\sss this map\dss is\dss a\sss homotopy\sss equivalence.\oss
Some weaker assumptions should\dss be sufficient\dss to ensure\sss that\trs
$X\qff \ttoo\qff A$\dss induces an isomorphism of\trs the fundamental\dss groups.\oss
In order\sss to find\sss $\pi_{\dff 1}\dff(\trf X\trf)$\sss one can\sss present\sss
$X$\sss as\sss $X\off =\off U\dff \cup\dff V$\dnsp,\oss where\sss
$U
\off =\off
C
\qff \cup\qff 
B_{\dff \bullet}\dff \times\dff [\trf 0\fff,\qff 2/3\trf)$
and\sss $V$\sss be\sss the image of\trs
$B_{\dff \bullet}\dff \times\dff (\trf 1/3\fff,\qff 1\trf]$\dss
in\dss $X$\nnsp.\oss
Then 
$U\dff \cap\dff V
\off =\off 
B_{\dff \bullet}\dff \times\dff (\trf 1/3\fff,\qff 2/3\trf)$\nnsp.\oss
In\sss the second special\sss case of\dss van\dss Kampen\dss
taking\sss as $U\fff,\qff V$\sss the images of\dss
$C_{\dff 1}\qff \cup\qff
B_{\dff 1}\dff \times\dff [\trf 0\fff,\qff 1\trf]\qff \cup\qff 
B_{\dff 2}\dff \times\dff (\trf 1/3\fff,\qff 1\trf]$\nnsp,\dss
$C_{\dff 2}\qff \cup\qff B_{\dff 2}\dff \times\dff [\trf 0\fff,\qff 2/3\trf)$
respectively\dss results in path-connected\sss intersection\sss $U\dff \cap\dff V$\dnsp.\oss

\myuppar{Non-abelian cohomology\sss in\sss dimensions $0\fff,\dff 1$\nnsp.}
P.\trs Olum\qss \cite{o}\qss introduced\qss \emph{singular\dss cohomology\sss groups}\dss
$H^{\dff 0}\dff(\trf X\fff,\qff A\dff;\qff \Pi\trf)$\nnsp,\dss 
$H^{\dff 1}\dff(\trf X\fff,\qff A\dff;\qff \Pi\trf)$\sss with possibly non-abelian\sss
groups of\dss coefficients $\Pi$\nnsp.\oss
His\sss theory\sss is\sss modeled on\dss the\dss Eilenberg--Steenrod\dss axiomatic approach\sss
to\sss the (co)homology\dss theory and\sss intended\sss for applications\sss to\sss
homotopy classification of\dss mappings.\oss
One of\dss his main\sss results was a\dss Mayer--Vietoris\dss sequence for such
cohomology\sss groups.\oss 
As an application of\dss this\dss
Mayer--Vietoris\dss sequence\dss Olum\dss presented a new proof\dss
of\dss the second,\oss Seifert--like,\oss special\sss case of\dss van Kampen\dss theorem.\oss
Adams\qss \cite{a}\qss called\dss this proof\qss \emph{``simple and conceptual''}.\oss

P.\dss Olum\dss provided\sss neither a new proof\dss of\sss the first\sss special\sss
case of\trs van Kampen\dss theorem,\oss nor even a new computation of\dss 
$\pi_{\dff 1}\dff(\trf S^{\fff 1}\trf)$\nnsp.\oss 
As evidenced\dss by\dss R.\dss Crowell\qss \cite{c}\qss and\dss R.\dss Fox\qss \cite{f},\oss
the computation of\sss $\pi_{\dff 1}\dff(\trf S^{\fff 1}\trf)$ wasn't\sss considered\sss
a worthwhile problem at\sss the\sss time,\oss in contrast\sss with\sss looking\sss
for conceptual\dss proofs of\dss Seifert\dss and\dss van\dss Kampen\dss theorems.\oss
Nevertheless,\oss seven\sss years\sss later\dss R.\dss Brown\dss opened\sss his paper\qss \cite{b1}\qss
with\sss the exclamation\qss
\emph{``We present\sss another proof\trs that
$\pi_{\dff 1}\dff(\trf S^{\fff 1}\trf)\off =\off \zzz$\nsp\textup{!}''}.\oss 
Actually,\oss R.\dss Brown\qss \cite{b1}\qss adapted\dss the method of\trs Olum\dss
to\sss find\sss $\pi_{\dff 1}\dff(\trf X\trf)$ when\sss $X\off =\off U\dff \cup\dff V$\dnsp,\oss
the interiors of\dss $U\fff,\qff V$\sss cover\sss $X$\nnsp,\oss
the subspaces $U\fff,\qff V$ are simply-connected,\oss
and\sss the intersection\sss $U\dff \cap\dff V$\sss consists of\dss $n\qff +\qff 1$\sss
path-components.\oss
Namely,\oss under\sss these assumptions $\pi_{\dff 1}\dff(\trf X\trf)$\sss
is\dss a free group on {\nsp}$n$ generators.\oss
The proof\trs is\dss again simple and conceptual\sss
and\sss includes an elegant\sss abstract\sss nonsense style argument\sss of\trs
Adams.\oss
We will\dss use\sss this argument\sss for\sss the same purpose.\oss
See\dss Theorems\qss \ref{pseudo-circle}\qss and\qss \ref{pseudo-wedge-of-circles}.\oss\vspace{-1pt}

R.\dss Brown\dss did not\sss develop\dss further these methods and embraced,\oss
starting with\dss the paper\qss \cite{b2}\qss and\sss the first\sss edition of\dss the book\qss
\cite{b3},\oss the ideology of\dss groupoids.\oss
Later on\sss he wrote about\qss \emph{``all\dss  the turgid stuff on nonabelian cohomology''.}\oss
See\qss \cite{bhs},\oss Section\qss 1.5.\oss
In\sss the present\sss paper we\sss follow\qss ``the road not\dss taken''\qss
by\dss R.\dss Brown and\sss use\sss the non-abelian cohomology.\vspace{-1pt}

\myuppar{Crowell{\fff}--{\fff}Fox\qss version of\qss
Seifert--van Kampen\qss theorem.}
The results of\trs Seifert\qss \cite{s}\qss
and\trs van Kampen\qss \cite{vk}\qss were stated\sss
in\sss terms of\dss generators and\sss relations.\oss
In early\dss 1950ies\dss Fox\dss reformulated\dss 
the standard\dss Seifert--van Kampen\dss theorem\qss 
(i.e.\qss the second special\sss case of\trs van Kampen)\qss 
in\sss terms of\dss direct\dss limits of\dss groups.\oss
The proof\dss was worked out\dss by\dss R.\dss Crowell\qss \cite{c}\qss
and\dss led\dss to a more general\dss result\sss about\sss unions
of\dss several\sss subsets.\oss
Suppose\sss that\sss a\sss topological\sss space $X$\sss is\dss presented\sss
as\sss the union\sss $X\off =\off \cup_{\dff i}\dff U_{\fff i}$\sss
of\dss open\sss path-connected subsets $U_{\fff i}\qff \subset\qff X$\nnsp.\oss
Suppose further\sss that\sss every subset\sss $U_{\fff i}$ 
contains a fixed\dss base point $b\qff \in\qff X$\nnsp,\oss
and\sss that\sss the family of\dss subsets $U_{\fff i}$\sss is\dss
closed under\sss finite intersections.\oss
Then\sss $\pi_{\dff 1}\dff(\trf X\fff,\qff b\trf)$\sss is\dss a direct\dss limit\sss
of\dss groups\sss $\pi_{\dff 1}\dff(\trf U_{\fff i}\fff,\qff b\trf)$\sss
and\sss homomorphisms\sss induced\dss by\sss inclusions of\dss the form\sss
$U_{\fff i}\qff \ttoo\qff U_{j}$\nsp.\oss

An examination of\dss the proof\dss shows\sss that\sss that\sss it\dss is\dss
sufficient\sss to assume\sss that\sss the family of\dss sets $U_{\fff i}$\sss
is\dss closed\sss under\sss taking\sss the intersection of\dss pairs of\dss sets
and\sss that\sss intersections of\dss $\leq\qff 4$\sss sets $U_{\fff i}$\sss
are path-connected.\oss
R.\dss Brown\dss and\dss A.R.\dss Salleh\qss \cite{bs},\oss
working\sss in\sss the groupoid\dss language,\oss showed\dss in\dss 1984\dss that\sss the\sss
last\sss condition can\sss be relaxed.\oss
Namely,\oss it\dss is\dss sufficient\sss to assume\sss that\sss
intersections of\dss $\leq\qff 3$\sss sets $U_{\fff i}$\sss are path-connected.\oss
The groupoid\dss language\dss is\dss irrelevant\sss for\sss this improvement.\oss
The corresponding argument\dss is\dss based on\sss the fact\sss that\sss the\dss
Lebesgue\dss covering dimension of\dss the disc\dss is\sss $\leq\qff 2$\nnsp.\oss
A proof\dss written\sss in\sss the usual\dss language of\dss groups\dss is\dss
contained\sss in\dss Hatcher's\trs textbook\qss \cite{h}.\oss
See\qss \cite{h},\oss Theorem\qss 1.20.\oss
As we will\sss see in\sss the proof\dss of\trs Theorem\qss \ref{olum-family},\oss
the\dss Lebesgue dimension of\dss the disc\dss is\dss also irrelevant.\oss\vspace{-1pt}

\myuppar{The present\sss paper\halfff.}
Our\sss first\sss goal\dss is\dss to present\sss a simple and elementary\sss proof\dss
of\dss the standard\dss Seifert--van Kampen\dss theorem\dss based on\sss the ideas of\trs
Olum\qss \cite{o}.\oss
In contrast\sss with\dss Olum,\oss we do not\sss discuss analogues of\dss the\dss
Eilenberg--Steenrod\dss axioms and\dss the\dss Mayer--Vietoris\dss sequence.\oss
Instead,\oss we work\sss mostly\sss with non-abelian cochains an cocycles
and use\sss the standard\dss tool\sss of\dss subdividing\sss the unit\sss square
into small\sss squares.\oss
The unpleasant\sss part\sss of\dss standard\dss proofs,\oss
the need\dss to keep\sss track of\dss paths connecting\sss the subdivision
points\sss to\sss a base point,\oss
is\dss absorbed\dss by\sss the\sss notion of\dss cohomologous cocycles.\oss
The resulting\sss proof\dss requires the same prerequisites as\sss 
the proof\dss in any\sss modern\sss textbook.\oss
It\sss occupies\dss Sections\qss \ref{cochains} -- \ref{svk-main}.\oss
The key\sss geometric result\dss is\dss Theorem\qss \ref{short-long},\oss
proved\dss by\sss subdividing\sss the unit\sss square.\oss

In\dss Section\qss \ref{svk-families}\qss we extend\dss the methods of\trs
Section\qss \ref{svk-main}\qss in order\sss to prove\sss the\dss
Crowell{\fff}--{\fff}Fox\qss version of\trs Seifert--van Kampen\qss theorem
and\dss its\dss Brown-Salleh\dss improvement.\oss
The proof\dss shows\sss that\sss the reason behind\dss the condition\sss that\sss
triple intersections are path-connected\dss is\dss purely combinatorial.\oss
There\dss is\dss no need\dss to subdivide\sss the unit\sss square into small\sss
rectangles without\sss fourfold\sss intersections.\oss
This contrasts with\qss \cite{bs},\qss 
and with\qss \cite{h},\oss the proof\dss of\qss Theorem\qss 1.20.\oss

Finally,\oss in\dss Section\qss \ref{relative-cohomology}\qss we extend\dss
the methods of\dss the previous sections\sss to\sss determine\sss
$\pi_{\dff 1}\dff(\trf X\trf)$\sss when
$X\off =\off U\dff \cup\dff V$\dnsp,\oss
the subsets $U\fff,\qff V$\sss are open,\oss
and\sss the intersection\sss $U\dff \cap\dff V$\sss consists of\dss finitely\sss many\sss
path-components.\oss
This proves\sss the open sets version of\dss the main\sss results of\trs van Kampen\qss
(and,\oss in\sss particular\halfff,\oss computes\sss the fundamental\sss group of\dss the circle).\oss
We deal\sss with\sss the case when\sss $U\dff \cap\dff V$\sss consists of\dss two
path-components\sss in\dss Theorem\qss \ref{open-van-kampen},\oss
and\sss with\sss the general\sss case in\dss Theorem\qss \ref{open-van-kampen-general}.\oss
The same methods work in\sss the case of\dss infinitely many components.\oss

\myuppar{Small\sss simplices and short\sss paths.}
Let\sss $\mathcal{U}$\sss be an open covering of\dss a\sss topological\sss space $X$\nnsp.\oss
Let\sss us call\sss a singular simplex\sss in $X$\sss \emph{small}\pss if\dss its
image\dss is\dss contained\sss in some $U\qff \in\qff \mathcal{U}$\dnsp.\oss
The usual\sss definition of\dss the singular homology and cohomology groups
with coefficients in an abelian group $\Pi$ can be modified\dss by 
considering only small\sss singular simplices.\oss
We will\dss indicate\sss this modification by\sss 
the subscript\sss $\mathcal{U}$\dnsp,\qss as\sss in\sss 
$H^{\fff n}_{\dff \mathcal{U}}\dff(\trf X\fff,\qff A\dff;\qff \Pi\trf)$\nnsp.\oss
There\dss is\dss an obvious map\vspace{3pt}
\begin{equation}
\label{small-cohomology}
\quad
H^{\fff n}\dff(\trf X\fff,\qff A\dff;\qff \Pi\trf)
\qff \ttoo\qff
H^{\fff n}_{\dff \mathcal{U}}\dff(\trf X\fff,\qff A\dff;\qff \Pi\trf)
\qff,
\end{equation}

\vspace{-12pt}\vspace{3pt}
and\dss by a fundamental\dss theorem,\oss
essentially due\sss to\dss Eilenberg\qss \cite{e},\oss
this map\sss is\sss an\sss isomorphism.\oss
See,\oss for example,\oss Proposition\qss 2.21\qss in\dss Hatcher's\dss textbook\qss \cite{h}.\oss
This results\sss is\dss the key\sss geometric step\sss to\sss the basic results
of\dss the singular (co)homology\sss theory such as\sss the excision property
and\sss the\dss Mayer--Vietoris\dss sequence.\oss
The geometric part\sss of\dss the proof\dss is\dss based on subdividing simplices into smaller ones,\qss
the barycentric subdivision\dss being\sss the standard\dss tool.\oss
The algebraic part\dss is\dss an elegant\sss construction of\dss chain homotopies
due\sss to\dss Eilenberg\qss \cite{e}.\oss

P.\dss Olum\dss observed\sss that\sss the same arguments work\sss for cohomology sets
with non-abelian coefficients $\Pi$ in dimensions $0\fff,\qff 1$\nnsp.\oss
In\sss fact,\oss the algebraic part\sss of\dss the proof\dss is\dss even simpler\halfff.\oss
See\qss \cite{o},\oss the proof\dss of\qss (2.5).\oss
As in\sss the (co)homology\sss theory,\oss
Olum\dss uses\sss the iterated\dss barycentric subdivisions of\dss simplices of\dss dimension $\leq\qff 2$\nnsp.\oss
The present\sss author believes\sss that\sss this unification with\sss the (co)homology\sss theory\sss 
is\dss an\sss important\sss advantage of\dss the non-abelian cohomology approach\sss 
to\dss Seifert--van Kampen\dss theorems.\oss

In\sss the context\sss of\dss fundamental\sss groups
the barycentric subdivisions are not\sss quite natural,\oss 
at\sss least\dss if\dss the\sss theory of\dss 
fundamental\sss groups precedes\sss the (co)homology\sss theory,\oss
as it\dss is\dss usually\sss the case.\oss 
By\sss this reason\sss we replaced\sss singular simplices by\sss paths
and\sss homotopies.\oss
We call\sss a path or a homotopy\qss \emph{short}\pss if\dss its image\dss 
is\dss contained\sss in some $U\qff \in\qff \mathcal{U}$\dnsp.\oss
The iterated\dss barycentric subdivision of\dss triangles are replaced\sss by\sss the much
simpler\sss tool\sss of\dss subdividing a square into smaller squares.\oss 
See\sss the proof\dss of\qss Theorem\qss \ref{short-long},\oss the analogue of\dss
the isomorphism\qss (\ref{small-cohomology}).\oss

\myuppar{Notations and conventions.}
The interval\dss $[\dff 0\fff,\qff 1\dff]$\dss is\dss often\sss denoted\dss by\dss $I$\nnsp.\oss
A\qss \emph{path}\qss in a\sss topological\sss spase\sss $X$\sss is\dss
a map\dss $I\qff \ttoo\qff X$\nnsp.\oss
If\dss $p\fff,\qff q$ are\sss two paths and\sss $p\dff(\dff 1\dff)\off =\off q\dff(\dff 0\dff)$\nnsp,\oss
then\sss the product\sss $p\dff \cdot\dff q$\sss is\dss defined\sss in\sss the usual\sss
manner\sss by\sss following first\sss $p$ and\sss then $q$\nnsp.\oss
If\dss $p$\sss is\dss a\sss path,\oss then\dss 
$\overline{p}$\qss
is\dss defined\dss by\trs 
$\overline{p}\trf(\dff s\trf)\off =\off p\dff(\dff 1\qff -\qff s\trf)$\nnsp.\oss
\emph{Homotopies of\qss paths}\qss are assumed\dss to be\qss 
\emph{homotopies relatively}\qss to\sss the boundary\sss 
$\{\trf 0\fff,\qff 1\trf\}$\nnsp.\oss
We write\dss $p\qff \sim\qff q$\dss 
when\sss paths\dss $p\fff,\pff q$\dss are homotopic.\oss

A\qss \emph{reparametrization}\qss of\dss a path\dss
$p$ is\dss a\sss path of\trs the form\dss
$p\dff \circ\dff \varphi$ where\dss
$\varphi\dff \colon\dff I\qff \ttoo\qff I$\dss
is\dss a\sss map such\dss that
$\varphi\dff(\dff 0\dff)\off =\off 0\fff,\off
\varphi\dff(\dff 1\dff)\off =\off 1$\nnsp.\oss
Clearly,\qss any\sss reparametrization of $p$\sss is\dss homotopic\sss to $p$\nnsp.\oss

Let\sss
$p_{\dff 1}\dff,\off p_{\dff 2}\dff,\off \ldots\dff,\off p_{\dff k}$\dss
be paths\sss in\sss $X$\sss such\dss that\dss the products
$p_{\dff i}\dff \cdot\dff p_{\dff i\dff +\dff 1}$\dss are defined,\oss
i.e.\qss such\sss that\sss
$p_{\dff i}\dff(\dff 1\dff)\off =\off p_{\dff i\dff +\dff 1}\dff(\dff 0\dff)$\sss
for every\dss $i\qff \leq\qff k\qff -\qff 1$\nnsp.\oss
Let\dss us define\dss 
$p\off =\off
p_{\dff 1}\dff \cdot\dff p_{\dff 2}\dff \cdot\dff \ldots\dff \cdot\dff p_{\dff k}$\dss
by\dss the rule\vspace{3pt}
\[
\quad
p\dff(\dff s\trf)
\off =\off
p_{\dff i}\dff\bigl(\trf k\dff s\qff -\qff i\qff +\qff 1\trf\bigr)\quad\
\mbox{for}\quad\
(\trf i\qff -\qff 1\trf)/k\qff \leq\qff s\qff \leq\qff i\fff/k
\pff.
\]

\vspace{-12pt}\vspace{3pt}
Clearly\halfff,\pss 
$p_{\dff 1}\dff \cdot\dff p_{\dff 2}\dff \cdot\dff \ldots\dff \cdot\dff p_{\dff k}$\dss 
differs\sss by\sss a\sss reparametrization\dss
from each\dss product\sss obtained\dss
by\dss placing\dss parentheses\sss into\sss the expression\dss
$p_{\dff 1}\dff \cdot\dff p_{\dff 2}\dff \cdot\dff \ldots\dff \cdot\dff p_{\dff k}$\nsp.\oss
Therefore\sss the product\sss
$p_{\dff 1}\dff \cdot\dff p_{\dff 2}\dff \cdot\dff \ldots\dff \cdot\dff p_{\dff k}$\dss
is\dss homotopic\sss to every\sss such\sss product\sss
and\sss may serve as a\sss partial\dss replacement\sss of\dss associativity\halfff.\oss

\newpage
\mysection{Group-valued\qss cochains\qss in\qss dimensions\dss $0$\dss and\dss $1$}{cochains}

\myuppar{Cochains and cocycles.}
Let\dss $X$\dss be a\sss path-connected\sss space,\qss $Y\qff \subset\qff X$\nnsp,\oss 
and $b\qff \in\qff Y$\nnsp.\oss
Let\sss us\sss fix a\sss group\dss $G$\sss
and denote by $1$\sss the unit\sss of\trs $G$\nnsp.\oss
Let\sss us denote by $P\dff(\trf X\trf)$\sss the set\sss of\dss all\dss paths in\dss $X$\nnsp,\oss
i.e.\qss the set\sss of\dss all\sss continuous maps\dss 
$[\dff 0\fff,\qff 1\dff]\qff \ttoo\qff X$\nnsp.\oss

A \emph{$0$\dnsp-cochain}\qss of\trs the pair\sss $(\trf X\fff,\pff Y\trf)$\dss  
is\dss a map\dss $c\dff \colon\dff X\qff \ttoo\qff G$\dss
such\dss that\trs $c\trf(\dff y\trf)\off =\off 1$\dss if\trs $y\qff \in\qff Y$\nnsp,\oss
and a \emph{$1$\dnsp-cochain}\qss of\dss $(\trf X\fff,\pff Y\trf)$\dss  
is\dss a map\dss $u\dff \colon\dff P\dff(\trf X\trf)\qff \ttoo\qff G$\dss
such\dss that\trs $u\trf(\dff p\trf)\off =\off 1$\dss 
if\trs $p\qff \in\qff P\dff(\trf Y\trf)$\nnsp.\oss
The set\sss of\dss all\dss $n$\dnsp-cochains of\dss $(\trf X\fff,\pff Y\trf)$\nnsp,\oss
where $n\off =\off 0$\sss or\sss $1$\nnsp,\oss is\dss denoted\dss by\trs
$C^{\dff n}\dff(\trf X\fff,\pff Y\trf)$\nnsp.\oss
Of\dss course,\pss $C^{\dff n}\dff(\trf X\fff,\pff Y\trf)$\dss
depends also on\dss $G$\nnsp.\oss
A $1$\dnsp-cochain $u$
is\dss called a \emph{cocycle}\pss 
if\sss 
$u\trf(\dff p\trf)\off =\off u\trf(\dff q\trf)$\dss
when\sss $p\off \sim\off q$\sss and 
$u\trf(\dff p\dff \cdot\dff q\trf)
\off =\off 
u\trf(\dff p\trf)\dff \cdot\dff u\trf(\dff q\trf)$\dss
when\dss the product\dss $p\dff \cdot\dff q$\dss is\dss defined.\oss
The set\sss of\dss all\sss cocycles of\dss $(\trf X\fff,\pff Y\trf)$\dss 
is\dss denoted\dss by\dss
$Z^{\dff 1}\dff(\trf X\fff,\pff Y\trf)$\nnsp.\oss

The case of\dss $Y\off =\off \{\dff b\trf\}$\sss is\dss the most\sss
important\sss one.\oss 
We will\sss abbreviate\sss
$(\trf X\fff,\pff \{\dff b\trf\}\trf)$\sss
as\sss $(\trf X\fff,\pff b\trf)$\nnsp.\oss
The case of\dss discrete subsets\sss $Y$\sss is\dss also important.\oss 
In\dss this case every\sss $p\qff \in\qff P\dff(\trf Y\trf)$\sss is\dss a constant\sss path
and\dss hence\sss $p\off \sim\off p\dff \cdot\dff p$\nnsp.\oss
Then\sss 
$u\trf(\dff p\trf)\off =\off 1$\dss 
for every $1$\dnsp-cocycle $u$ and\sss $p\qff \in\qff P\dff(\trf Y\trf)$\nnsp.\oss

\myuppar{Cohomology\halfff.}
The point-wise multiplication of\dss maps $X\dff \ttoo\dff G$
turns $C^{\dff 0}\dff(\trf X\fff,\qff Y\trf)$ into a\sss group.\oss
The group\dss
$C^{\dff 0}\dff(\trf X\fff,\pff Y\trf)$ acts on
$C^{\dff 1}\dff(\trf X\fff,\pff Y\trf)$\sss by\dss the rule\vspace{0pt}
\[
\quad
(\trf c\dff \bullet\dff u\trf)\trf(\dff p\trf)
\off =\off
c\trf(\trf p\dff(\dff 0\dff)\trf)\dff \cdot\dff
u\trf(\dff p\trf)\dff \cdot\dff
c\trf(\trf p\dff(\dff 1\dff)\trf)^{\dff -\dff 1}
\pff,
\]

\vspace{-12pt}
where\dss
$c\qff \in\qff C^{\dff 0}\dff(\trf X\fff,\pff Y\trf)$\nnsp,\qss
$u\qff \in\qff C^{\dff 1}\dff(\trf X\fff,\pff Y\trf)$\nnsp,\oss
and\dss $p\qff \in\qff P\dff(\trf X\trf)$\nnsp.\oss
An\sss immediate verification shows\sss that\dss
$c\dff \bullet\dff (\trf d\dff \bullet\dff u\trf)
\off =\off
(\dff c\dff \cdot\dff d\trf)\dff \bullet\dff u$\sss
and\dss hence\dss
$(\dff c\fff,\pff u\trf)\qff \longmapsto\qff c\dff \bullet\dff u$\dss
is\dss indeed an action.\oss
Another\sss easy verification shows\sss that\dss
$c\dff \bullet\dff u$\sss is\dss a cocycle\sss
if\sss $u$\sss is\dss a cocycle.\oss
Two cocycles\sss
$z\fff,\pff u$\sss
are called\qss \emph{cohomologous}\pss if\qss
they\sss belong\dss to\sss the same orbit\sss of\trs the action of\dss
$C^{\dff 0}\dff(\trf X\fff,\pff Y\trf)$\nnsp,\oss
i.e.\qss if\qss $z\off =\off c\dff \bullet\dff u$\qss
for some\sss $c\qff \in\qff C^{\dff 0}\dff(\trf X\fff,\pff Y\trf)$\nnsp.\qss
Being cohomologous\dss is\dss an equivalence relation on\dss the set\dss
$Z^{\dff 1}\dff(\trf X\fff,\pff Y\trf)$\nnsp.\oss
The set\sss of\dss equivalence classes\dss is\dss denoted\sss by\sss
$H^{\dff 1}\dff(\trf X\fff,\pff Y\trf)$\nnsp.\oss
The equivalence class of\dss a cocycle $u$\sss is\dss called\sss
the\qss \emph{cohomology\sss class}\qss of\sss $u$
and denoted\dss by $\hclass{u\fff}$\nnsp.\oss
If\sss $G$\sss is\dss not\sss abelian,\dss
$H^{\dff 1}\dff(\trf X\fff,\pff Y\trf)$ 
has no\sss natural\sss group structure,\qss
but\dss has a distinguished element,\qss
called\qss \emph{trivial}\pss cohomology\sss class,\oss namely\halfff,\oss
the equivalence class $1$ of\dss the cocycle $\mathbb{1}$ 
mapping every\sss path\dss to $1\qff \in\qff G$\nnsp.\oss

\myuppar{Cohomology\sss classes and\dss homomorphisms.}
It\dss turns out\dss that\dss 
there\dss is\dss a canonical\dss bijection\dss between\dss
$H^{\dff 1}\dff(\trf X\fff,\pff b\dff)$\dss
and\dss the set\trs
$\hom (\trf \pi_{\dff 1}\dff(\trf X\fff,\qff b\trf)\fff,\pff G\trf)$\dss
of\trs homomorphisms\dss
$\pi_{\dff 1}\dff(\trf X\fff,\qff b\trf)
\qff \ttoo\qff G$\nnsp.\oss
Let\trs $L\trf(\trf X\fff,\qff b\trf)$\sss be\sss 
the set\sss of\trs loops in\dss $X$\dss based at\sss $b$\nnsp,\oss
and\dss let\dss us\dss interpret\dss homomorphisms\dss
$\pi_{\dff 1}\dff(\trf X\fff,\qff b\trf)
\qff \ttoo\qff G$\dss
as maps\sss 
$h\dff \colon\dff
L\trf(\trf X\fff,\qff b\trf)
\qff \ttoo\qff G$\dss
such\dss that\dss
$h\trf(\dff p\trf)\off =\off h\trf(\dff q\trf)$\dss
if\sss $p\off \sim\off q$\nnsp,\oss 
and such\dss that\sss 
$h\trf(\dff p\dff \cdot\dff q\trf)
\off =\off 
h\trf(\dff p\trf)\dff \cdot\dff h\trf(\dff q\trf)$\sss 
for every\sss $p\fff,\pff q\qff \in\pff L\trf(\trf X\fff,\qff b\trf)$\nnsp.\oss
Let\dss $u\qff \in\qff Z^{\dff 1}\dff(\trf X\fff,\pff b\dff)$\nnsp.\oss 
Since\sss $u$\sss is\dss a cocycle,\oss
the restriction of\dss $u$\sss to\dss 
$L\trf(\trf X\fff,\qff b\trf)\qff \subset\qff P\dff(\trf X\trf)$\dss
satisfies\sss the above conditions and\dss hence
defines a homomorphism\dss
$\rho\dff(\dff u\dff)\dff \colon\dff
\pi_{\dff 1}\dff(\trf X\fff,\qff b\trf)
\qff \ttoo\qff G$\nnsp.\oss
Since\sss elements\sss of\trs
$C^{\dff 0}\dff(\trf X\fff,\pff b\dff)$\dss
are\dss equal\dss to $1$ at\sss $b$\nnsp,\qss
$\rho\dff(\dff u\dff)$ depends only\sss on\dss the cohomology\sss class $\hclass{u\fff}$\nnsp,\oss
and we get\sss a map\vspace{0.75pt}
\[
\quad
\rho\qff \colon\qff
H^{\dff 1}\dff(\trf X\fff,\pff b\dff)
\off \ttoo\off
\hom (\trf \pi_{\dff 1}\dff(\trf X\fff,\qff b\trf)\fff,\pff G\trf)
\pff.
\]

\vspace{-12pt}\vspace{0.75pt}
In order\dss to construct\dss a map in\dss the opposite direction,\oss
let\dss us\sss choose for every\sss point\sss $x\qff \in\qff X$\sss
a path $s_{\fff x}$ such\dss that\sss
$s_{\fff x}\trf(\dff 0\dff)\off =\off b$
and\sss
$s_{\fff x}\trf(\dff 1\dff)\off =\off x$\nnsp,\oss
and choose as $s_{\dff b}$\sss
the constant\sss path\sss with\sss the value $b$\nnsp.\oss
Let\dss us\dss define a map\dss
$l\dff \colon\dff
P\dff(\trf X\trf)\qff \ttoo\qff L\trf(\trf X\fff,\qff b\trf)$\dss
by\dss the rule\sss \vspace{3pt}\vspace{-0.375pt}
\[
\quad
l\trf(\dff p\trf)
\off =\off 
s_{\dff p\dff(\dff 0\dff)}
\dff \cdot\dff 
p
\dff \cdot\pff 
\overline{s_{\dff p\dff(\dff 1\dff)}}
\off.
\]

\vspace{-12pt}\vspace{3pt}\vspace{-0.375pt}
If\sss $p\qff \sim\qff q$\nnsp,\oss then,\oss in\dss particular\halfff,\qss
$p\dff(\dff 0\dff)\off =\off q\dff(\dff 0\dff)$
and\sss
$p\dff(\dff 1\dff)\off =\off q\dff(\dff 1\dff)$
and\dss hence\sss
$l\trf(\dff p\trf)\off \sim\off l\trf(\dff q\trf)$\nnsp.\oss
If\sss $p$\sss is\sss a constant\sss path,\oss
then\sss $p\dff(\dff 0\dff)\off =\off p\dff(\dff 1\dff)$\sss
and\sss
$l\trf(\dff p\trf)$\sss is\dss homotopic\sss to\sss the constant\sss loop.\oss
If\dss the product\sss $p\dff \cdot\dff q$\sss is\dss defined,\oss
then $p\dff(\dff 1\dff)\off =\off q\dff(\dff 0\dff)$
and\dss hence\dss
$\overline{s_{\dff p\dff(\dff 1\dff)}}
\qff \cdot\dff
s_{\dff q\dff(\dff 0\dff)}
\off =\off
\overline{s_{\dff p\dff(\dff 1\dff)}}
\qff \cdot\dff
s_{\dff p\dff(\dff 1\dff)}$\dss
is\dss homotopic\sss to a constant\dss path.\oss
It\dss follows\dss that\sss
$l\trf(\dff p\dff \cdot\dff q\trf)
\off \sim\off 
l\trf(\dff p\trf)\dff \cdot\dff
l\trf(\dff q\trf)$\nnsp.\oss

Let\sss
$h\dff \colon\dff
L\trf(\trf X\fff,\qff b\trf)
\qff \ttoo\qff G$\sss
be\dss a map satisfying\dss the above conditions.\oss
Then\dss the properties of\sss $l$\sss
imply\dss that\dss the map\sss
$h^{\dff \sim}\dff \colon\dff
P\dff(\trf X\trf)\qff \ttoo\qff G$\sss 
defined\dss by\sss 
$h^{\dff \sim}\dff(\dff p\trf)
\off =\off
h\trf(\trf l\trf(\dff p\trf)\trf)$\sss
is\dss a cocycle.\oss
Taking\dss the cohomology\sss classes of\dss the cocycles\dss $h^{\dff \sim}$\sss
leads\sss to a map\vspace{3pt}
\[
\quad
\varepsilon\qff \colon\qff
\hom (\trf \pi_{\dff 1}\dff(\trf X\fff,\qff b\trf)\fff,\pff G\trf)
\off \ttoo\off
H^{\dff 1}\dff(\trf X\fff,\pff b\dff)
\pff.
\]

\vspace{-12pt}\vspace{3pt}
\mypar{Lemma.}{hom-cohomology}
\emph{Both maps\dss $\rho$\dss and\dss $\varepsilon$\dss
are bijections and\dss $\varepsilon\off =\off \rho^{\dff -\dff 1}$\dnsp.\oss}

\proof
Since\sss $s_{\dff b}$\sss is\dss the constant\sss path,\pss
$l\trf(\dff p\trf)\off =\off
s_{\dff b}
\dff \cdot\dff 
p
\dff \cdot\qff 
\overline{s_{\dff b}}
\off\off \sim\off\qff
p$\dss
for every\dss
$p\qff \in\qff L\trf(\trf X\fff,\qff b\trf)$\dss
and\dss hence\sss the restriction of\dss 
$h^{\dff \sim}$\dss 
to\dss $L\trf(\trf X\fff,\qff b\trf)$\dss
is\dss equal\dss to\sss $h$\nnsp.\oss
Therefore\dss $\rho\dff \circ\dff \varepsilon$\dss is\dss the identity\sss map.\oss
It\dss remains\sss to show\dss that\trs 
$\varepsilon\dff \circ\dff \rho$\dss 
is\dss the identity\sss map.\oss
Let\dss $u\qff \in\qff Z^{\dff 1}\dff(\trf X\fff,\pff b\dff)$\nnsp.\oss
Let\dss us consider\dss its restriction\dss $h$\dss to\dss
$L\trf(\trf X\fff,\qff b\trf)$\dss and\dss the corresponding cocycle\dss $h^{\dff \sim}$\dnsp.\oss
If\trs $p\qff \in\qff P\dff(\trf X\trf)$\nnsp,\oss
then\vspace{3pt}
\[
\quad
h^{\dff \sim}\dff(\dff p\trf)
\off =\off
h\trf(\trf l\trf(\dff p\trf)\trf)
\off =\off
u\trf(\trf l\trf(\dff p\trf)\trf)
\]

\vspace{-39pt}
\[
\quad
\phantom{h^{\dff \sim}\dff(\dff p\trf)
\off =\off
h\trf(\trf l\trf(\dff p\trf)\trf)
\off }
=\off
u\trf\left(\qff 
s_{\dff p\dff(\dff 0\dff)}
\dff \cdot\dff 
p
\dff \cdot\pff 
\overline{s_{\dff p\dff(\dff 1\dff)}}
\pff\right)
\]

\vspace{-36pt}
\[
\quad
\phantom{h^{\dff \sim}\dff(\dff p\trf)
\off =\off
h\trf(\trf l\trf(\dff p\trf)\trf)
\off }
=\off
u\trf\left(\dff 
s_{\dff p\dff(\dff 0\dff)}
\trf\right)
\dff \cdot\dff 
u\trf\left(\dff p\trf\right)
\dff \cdot\pff 
u\trf\left(\pff
\overline{s_{\dff p\dff(\dff 1\dff)}}
\pff\right)
\]

\vspace{-36pt}
\[
\quad
\phantom{h^{\dff \sim}\dff(\dff p\trf)
\off =\off
h\trf(\trf l\trf(\dff p\trf)\trf)
\off }
=\off
u\trf\left(\dff 
s_{\dff p\dff(\dff 0\dff)}
\trf\right)
\dff \cdot\dff 
u\trf\left(\dff p\trf\right)
\dff \cdot\pff 
u\trf\left(\pff
s_{\dff p\dff(\dff 1\dff)}
\pff\right)^{\dff -\dff 1}
\pff.
\]

\vspace{-12pt}\vspace{3pt}
It\dss follows\dss that\trs
$h^{\dff \sim}\off =\off c\dff \bullet\dff u$\nnsp,\oss
where\dss $c$\dss is\dss the $0$\dnsp-cochain\dss
$x\off \longmapsto\off 
u\trf\left(\dff 
s_{\dff x}
\trf\right)$\nnsp.\oss
Since\sss $s_{\dff b}$\dss is\dss the constant\dss path with\dss
the value\sss $b$\sss and\dss
$u\qff \in\qff Z^{\dff 1}\dff(\trf X\fff,\pff b\dff)$\nnsp,\oss
the $0$\dnsp-cochain $c$\sss belongs\sss to\dss
$C^{\dff 0}\dff(\trf X\fff,\pff b\dff)$\nnsp.\oss
Therefore\sss the cohomology\sss class of\trs $h^{\dff \sim}$\dss
is\dss equal\dss to\sss the cohomology\sss class of\trs $u$\nnsp.\oss
Since $u$\sss was an arbitrary\sss cocycle,\oss
this implies\sss that\trs $\varepsilon\dff \circ\dff \rho$\dss 
is\dss the identity\sss map.\oss  \eproof

\mypar{Theorem.}{independence-of-paths}
\emph{The maps\trs $\rho$\sss and\trs $\varepsilon$\sss are
canonical\dss bijections\dss  
between\dss the sets\dss
$H^{\dff 1}\dff(\trf X\fff,\pff b\dff)$\dss and\dss 
$\hom (\trf \pi_{\dff 1}\dff(\trf X\fff,\qff b\trf)\fff,\pff G\trf)$\nnsp.\oss
In\dss particular\halfff,\pss $\varepsilon$\dss
does not\dss depend\dss
on\dss the choice of\trs paths\dss $s_{\dff x}$\nsp.\oss}

\proof
It\dss remains\sss only\dss to check\dss that\sss $\varepsilon$\sss
does not\dss depend on\dss the choice of\trs paths\dss $s_{\dff x}$\nsp.\oss
This follows from\dss the fact\dss that\sss definition of\dss 
$\rho$\dss does not\dss involve any\sss choices.\oss  \eproof

\mysection{Short\qss cochains}{short}

\myuppar{Short\dss paths and\sss short\dss homotopies.}
Let\dss us\dss fix an open covering $\mathcal{U}$ of\sss $X$\nnsp.\oss
Let\sss call\sss a path $p$ in $X$\dss 
\emph{short}\pss if\dss $p\dff(\trf I\trf)\qff \subset\qff U$\sss
for some $U\qff \in\qff \mathcal{U}$\dnsp.\oss
A\sss homotopy\sss 
$I\dff \times\dff [\dff 0\fff,\qff 1\dff]\qff \ttoo\qff X$\sss
is\dss called\qss \emph{short}\pss if\dss its image\dss is\dss contained\sss in\sss $U$\sss 
for some $U\qff \in\qff \mathcal{U}$\dnsp.\oss
Obviously\halfff,\oss a short\dss homotopy\dss is\dss a\sss homotopy\dss
between\sss two short\dss paths.\oss
One can\sss replace in\dss the definitions of\dss cochains and cocycles
arbitrary\dss paths by\sss short\dss paths and\sss 
arbitrary\dss homotopies\sss by\sss short\dss homotopies.\oss

In more details,\oss
let\sss $P_{\dff \mathcal{U}}\dff(\trf X\trf)$ 
be\sss the set\sss of\dss short\dss paths in $X$\nnsp.\oss
Every $0$\dnsp-cochain\dss should\dss be considered as\qss \emph{short}\qss
because every\sss point\dss belongs\sss to some
$U\qff \in\qff \mathcal{U}$\nnsp.\oss
A\qss \emph{short\sss $1$\dnsp-cochain}\qss 
of\dss the pair $(\trf X\fff,\pff Y\trf)$\sss 
is\dss a map\sss 
$u\dff \colon\dff P_{\dff \mathcal{U}}\dff(\trf X\trf)\qff \ttoo\qff G$\sss
such\dss that\sss $u\trf(\dff p\trf)\off =\off 1$\sss
if\sss $p\qff \in\qff P\dff(\trf Y\trf)$\nnsp.\oss
The set\sss of\dss all\sss short\sss $n$\dnsp-cochains,\oss
where $n\off =\off 0$\sss or\sss $1$\nnsp,\oss is\dss denoted\dss by\sss
$C^{\dff n}_{\dff \mathcal{U}}\dff(\trf X\fff,\pff Y\trf)$\nnsp.\oss
A short\sss $1$\dnsp-cochain $u$
is\dss called a \emph{short\sss cocycle}\pss 
if\sss 
$u\trf(\dff p\trf)\off =\off u\trf(\dff q\trf)$\dss
when\dss there exists a short\dss homotopy\dss
relatively\dss to\dss $\{\qff 0\fff,\qff 1\trf\}$\dss 
between $p$ and $q$\nnsp,\oss and\dss also\sss 
$u\trf(\dff p\dff \cdot\dff q\trf)
\off =\off 
u\trf(\dff p\trf)\dff \cdot\dff u\trf(\dff q\trf)$\dss
when\dss the product\dss $p\dff \cdot\dff q$\sss is\dss defined and\dss
the path $p\dff \cdot\dff q$\sss is\dss short.\oss
The set\sss of\dss all\sss short\sss cocycles\dss is\dss denoted\dss by\dss
$Z^{\dff 1}_{\dff \mathcal{U}}\dff(\trf X\fff,\pff Y\trf)$\nnsp.\oss

The action of\sss  $0$\dnsp-cochains
on\sss short $1$\dnsp-cochains\dss
is\dss defined exactly\sss as before,\oss
as also\sss the relation of\dss being\qss \emph{cohomologous}\qss cocycles.\oss
Being cohomologous\dss is\dss an equivalence relation,\oss and\dss
the set\sss of\dss equivalence classes of\sss short $1$\dnsp-cocycles\dss
is\dss denoted\dss by\sss
$H^{\dff 1}_{\dff \mathcal{U}}\dff(\trf X\fff,\pff Y\trf)$\nnsp.\oss
Restricting\sss maps
$P\dff(\trf X\trf)\qff \ttoo\qff G$
to $P_{\dff \mathcal{U}}\dff(\trf X\trf)$\dss 
leads\sss to maps\dss 
$F^{\dff 1}\dff(\trf X\fff,\pff Y\trf)
\off \ttoo\off
F^{\dff 1}_{\dff \mathcal{U}}\dff(\trf X\fff,\pff Y\trf)$\nnsp,\oss
where\dss $F\off =\off C$\nnsp,\dss $Z$\nnsp,\qss or\sss $H$\nnsp.

\mypar{Theorem.}{short-long}
\emph{For\qss $F\off =\off Z$\qss or\qss $H$\dss the map\qss
$F^{\dff 1}\dff(\trf X\fff,\pff Y\trf)
\off \ttoo\off
F^{\dff 1}_{\dff \mathcal{U}}\dff(\trf X\fff,\pff Y\trf)$\dss
is\dss a\dss bijection.\oss}

\proof
Let\dss us\dss deal\dss with cocycles first\halfff.\oss
Given a short\sss cocycle\sss $u$\nnsp,\oss let\dss us define a cochain\dss
$u^{\dff \sharp}$\dss as follows.\oss
Let\dss $p\dff \colon\dff I\qff \ttoo\qff X$\dss be a path.\oss
By\dss Lebesgue\qss lemma\dss there exists
a subdivision of\trs the interval $I$ into subintervals such\sss that
$p$ maps each of\dss these subintervals into some $U\qff \in\qff \mathcal{U}$\nnsp.\oss
In other\dss words,\oss there are numbers\sss 
$0
\off =\off
a_{\dff 1}\off <\off
a_{\dff 2}\off <\off
\ldots\off <\off
a_{\dff n}
\off =\off
1$\sss
such\dss that\sss 
$p\trf\trf [\dff a_{\dff i\dff -\dff 1}\dff,\pff a_{\dff i}\dff]\trf)$\sss
is\dss contained\sss in\dss some $U_{\dff i}\qff \in\qff \mathcal{U}$
for each\sss
$i\off =\off 1\fff,\pff 2\fff,\pff \ldots\fff,\pff n$\nnsp.\oss
For each such number\sss
$i$\sss
let\sss $p_{\dff i}\off =\off p\dff \circ\dff \varphi_{\dff i}$\nsp,\oss
where\dss 
$\varphi_{\dff i}\dff \colon\dff 
[\dff 0\fff,\qff 1\dff]
\qff \ttoo\qff
[\trf a_{\dff i\dff -\dff 1}\dff,\pff a_{\dff i}\qff]$\sss
is\dss an\sss increasing\sss homeomorphism.\oss
Then each\sss $p_{\dff i}$\sss is\dss a\sss short\dss path and\dss
hence\dss $u\trf(\dff p_{\dff i}\trf)$\dss is\dss defined.\oss
Let\vspace{3pt}
\begin{equation}
\label{u-sharp}
\quad
u^{\dff \sharp}\dff(\dff p\trf)
\off =\off
u\trf(\dff p_{\dff 1}\trf)\dff \cdot\dff
u\trf(\dff p_{\dff 2}\trf)\dff \cdot\dff
\ldots\dff \cdot\dff
u\trf(\dff p_{\dff n}\trf)
\pff.
\end{equation}

\vspace{-12pt}\vspace{3pt}
Let\dss us\sss check\dss that\trs $u^{\dff \sharp}\dff(\dff p\trf)$\dss
is\sss independent\dss from\dss the choices involved.\oss
To begin\dss with,\oss for each\sss $i$\sss 
different\sss choices of\trs $\varphi_{\dff i}$\sss
lead\dss to paths differing\dss by\sss a\sss reparametrization.\oss
Since\sss $p_{\dff i}$\sss is\dss a\sss short\dss path,\oss
a reparametrization of\dss $p_{\dff i}$\sss is\dss homotopic\sss to\sss
$p_{\dff i}$\sss by\sss a\sss short\dss homotopy\halfff.\oss
Since $u$ is\dss a\sss short $1$\dnsp-cocycle,\oss
this implies\sss that
$u\trf(\dff p_{\dff i}\trf)$ does not\sss depend on\dss
the choice of\dss $\varphi_{\dff i}$\nsp.\oss
It\dss follows\dss that $u^{\dff \sharp}\dff(\dff p\trf)$
does not\sss depend on\dss the choice of\trs homeomorphisms\sss $\varphi_{\dff i}$\nnsp.\oss
Replacing\sss  
$[\trf a_{\dff i\dff -\dff 1}\dff,\pff a_{\dff i}\qff]$\sss
by\sss its subdivisions\sss will\dss not\dss change\dss
$u^{\dff \sharp}\dff(\dff p\trf)$\dss because $u$\sss is\dss a short\sss cocycle.\oss
Since every\dss two subdivisions of\trs $I\off =\off [\dff 0\fff,\qff 1\dff]$\dss
have a common subdivision,\oss
it\dss follows\dss that\dss $u^{\dff \sharp}\dff(\dff p\trf)$\dss
does not\dss depend on\dss the subdivision used.\oss
Therefore\sss $u^{\dff \sharp}$\dss is\dss correctly\sss defined.\oss

Let\dss us\dss prove\sss that\dss $u^{\dff \sharp}$\dss 
is\dss a cocycle of\dss $(\trf X\fff,\pff Y\trf)$\nnsp.\oss
Clearly\halfff,\pss 
$u^{\dff \sharp}\dff(\dff p\trf)\off =\off 1$\dss if\trs
$p\qff \in\qff P\dff(\trf Y\trf)$\nnsp.\oss
The fact\dss that\sss
$u^{\dff \sharp}\dff(\dff p\dff \cdot\dff q\trf)
\off =\off 
u^{\dff \sharp}\dff(\dff p\trf)\dff \cdot\dff u^{\dff \sharp}\dff(\dff q\trf)$\dss
when\sss $p\dff \cdot\dff q$\dss is\dss defined\dss
follows\dss from\dss the independence of\dss
$u^{\dff \sharp}\dff(\dff p\dff \cdot\dff q\trf)$\dss on\sss the subdivision used.\oss
Let\dss us\dss
prove\sss that\sss
$u^{\dff \sharp}\dff(\dff p\trf)
\qff =\qff 
u^{\dff \sharp}\dff(\dff q\trf)$
when\trs $p\off \sim\off q$\nnsp.\oss
Let
$h\dff \colon\dff
I\dff \times\dff [\dff 0\fff,\qff 1\dff]
\qff \ttoo\qff X$\dss
be a\sss homotopy\dss between\dss 
$p$ and $q$ such\sss that\sss
$h\trf(\dff 0\fff,\qff t\trf)
\off =\off 
p\dff(\dff 0\trf)
\off =\off 
q\dff(\dff 0\trf)$ 
and\dss
$h\trf(\dff 1\fff,\qff t\trf)
\off =\off 
p\dff(\dff 1\trf)\off =\off q\dff(\dff 1\trf)$\dss
for every\sss $t\qff \in\qff [\dff 0\fff,\qff 1\dff]$\nnsp.

By\qss Lebesgue\qss lemma\dss there exists
a natural\sss number $n$ such\dss that\sss $h$\sss maps every\sss square\vspace{3pt}
\[
\quad
\left[\qff
\frac{\dff i\dff}{n}\qff,\off
\frac{i\qff +\qff 1}{n}
\qff\right]
\qff \times\qff
\left[\qff
\frac{\dff k\dff}{n}\qff,\off
\frac{k\qff +\qff 1}{n}
\qff\right]
\pff
\]

\vspace{-9pt}
with\dss
$i\fff,\pff k\off =\off 0\fff,\pff 1\fff,\pff \ldots\fff,\pff n\qff -\qff 1$\dss
into some\dss $U\qff \in\qff \mathcal{U}$\nnsp.\oss
Let\dss us consider segments of\trs the form\vspace{3pt}
\[
\quad
K
\off =\off
\left[\qff
\frac{\dff k\dff}{n}\qff,\off
\frac{k\qff +\qff 1}{n}
\qff\right]
\qff \times\qff
\frac{\dff i\dff}{n}
\quad\
\mbox{or}\quad\
\frac{\dff i\dff}{n}
\qff \times\qff
\left[\qff
\frac{\dff k\dff}{n}\qff,\off
\frac{k\qff +\qff 1}{n}
\qff\right]
\]

\vspace{-9pt}
with\dss 
$0\qff \leq\qff i\qff \leq\qff n$\dss 
and\dss
$0\qff \leq\qff k\qff \leq\qff n\qff -\qff 1$\nnsp.\oss
These segments are nothing else but\dss the sides of\trs the
above squares.\oss
For each such segment\dss $K$\dss let\trs
$\varphi_{\trf K}\dff \colon\dff
I\qff \ttoo\qss I\dff \times\dff I$\dss be\sss the path
defined\dss by\vspace{3pt}
\[
\quad
\varphi_{\trf K}\dff(\dff s\trf)
\off =\off
\left(\qff
\frac{\dff i\dff}{n}\dff,\off
\frac{k\qff +\qff s}{n}
\qff\right)
\quad\
\mbox{or}\quad\
\left(\qff
\frac{k\qff +\qff s}{n}\dff,\off
\frac{\dff i\dff}{n}
\qff\right)
\]

\vspace{-9pt}
respectively\halfff.\oss
Then\dss $\varphi_{\trf K}\dff(\trf I\trf)\off =\off K$\dss
and\dss $h_{\trf K}\off =\off h\dff \circ\trf \varphi_{\dff K}$\dss
is\dss a\sss short\dss path.\oss
If\trs $K$\dss is\dss contained\sss in\sss either\dss 
$0\dff \times\dff [\dff 0\fff,\qff 1\dff]$\dss
or\dss
$1\dff \times\dff [\dff 0\fff,\qff 1\dff]$\nnsp,\oss
then\dss $h_{\trf K}$\dss is\dss a\sss constant\dss path.\oss

\vspace{3pt}
\begin{figure}[h!]
\hspace*{10em}
\includegraphics[width=0.48\textwidth]{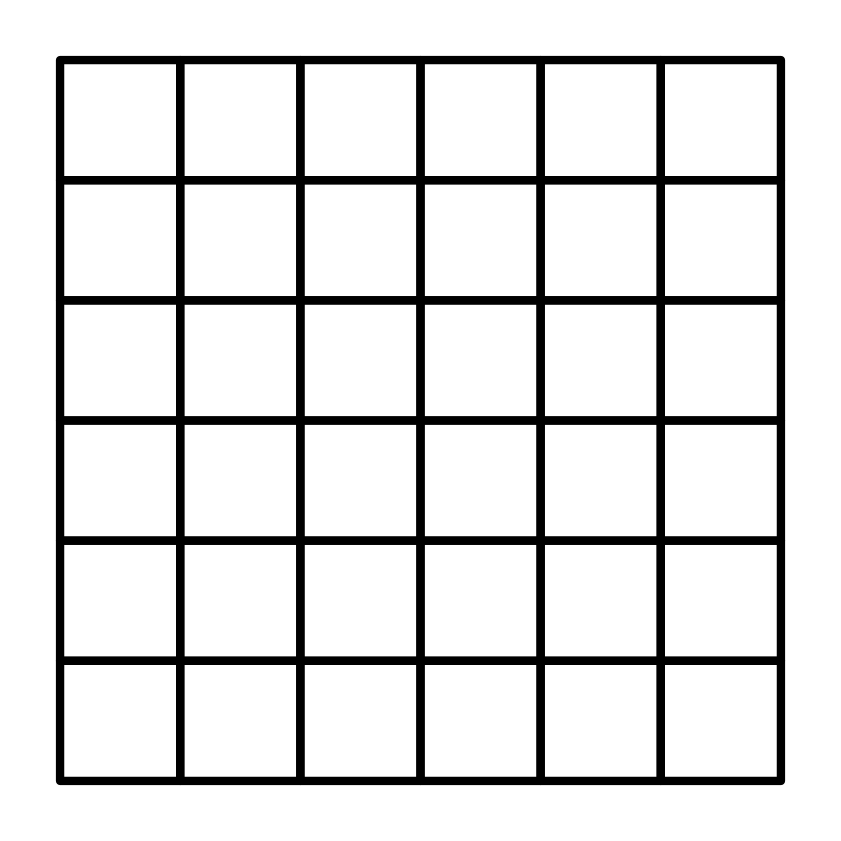}
\end{figure}\vspace{3pt}

Let\dss $p_{\fff j}
$\nnsp,\oss
where
$j\off =\off 0\fff,\pff 1\fff,\pff \ldots\fff,\pff n$\nnsp,\oss
be\sss the path in\dss $X$\dss defined\dss by\dss
$p_{\fff j}\dff(\dff s\trf)\off =\off h\trf(\dff s\fff,\pff j/n\trf)$\nnsp.\oss
Then\dss $p_{\dff 0}\off =\off p$\nnsp,\pss
$p_{\dff n}\off =\off q$\nnsp,\oss
and\dss hence\sss it\dss is\dss sufficient\dss to prove\sss that\dss 
$u^{\dff \sharp}\dff(\dff p_{\fff j}\trf)
\off =\off 
u^{\dff \sharp}\dff(\dff p_{\fff j\dff -\dff 1}\trf)$\dss
for every\sss
$j\off =\off 1\fff,\pff 2\fff,\pff \ldots\fff,\pff n$\nnsp.\oss
Let\sss us consider\sss the rectangle\dss
$Q_{\trf j}
\off =\off
I\dff \times\dff [\dff (\dff j\qff -\qff 1\dff)/n\fff,\pff j/n\dff]$\nnsp.\oss

\vspace{-6pt}
\begin{figure}[h!]
\hspace*{10em}
\includegraphics[width=0.48\textwidth]{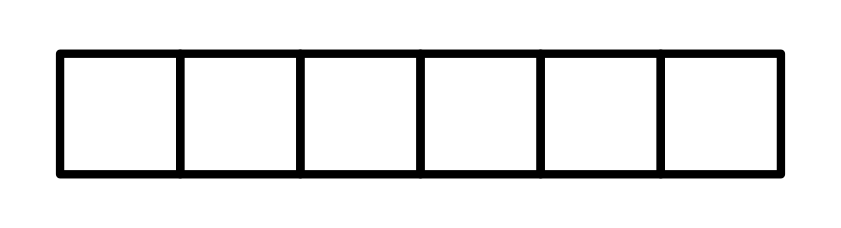}
\end{figure}
\vspace{-6pt}

Let\dss
$K_{\dff 1}\dff,\off K_{\trf 2}\dff,\off \ldots\dff,\off K_{\dff n}$\dss
be\sss the segments of\trs the above form contained\dss in\dss
$I\dff \times\dff j/n$\dss and\dss listed\dss from\dss left\dss to\sss right\halfff,\oss
and\dss let\dss 
$k_{\dff i}\off =\off h_{\trf K_{\dff i}}$.\oss
Similarly\halfff,\oss let\dss
$L_{\trf 1}\dff,\off L_{\qff 2}\dff,\off \ldots\dff,\off L_{\dff n}$\dss
be\sss the segments of\trs the above form contained\dss in\dss
$I\dff \times\dff (\dff j\qff -\qff 1\dff)\dff/n$\dss 
and\dss listed\dss from\dss left\dss to\sss right\halfff,\oss
and\dss let\dss 
$l_{\dff i}\off =\off h_{\trf L_{\trf i}}$.\oss
Then every\dss $k_{\dff i}$\dss and every\sss $l_{\dff i}$\sss 
is\dss a\sss short\dss path and\vspace{3pt}
\[
\quad
p_{\fff j}
\off =\off
k_{\dff 1}\dff \cdot\qff
k_{\qff 2}\dff \cdot\qff
\ldots\qff \cdot\qff
k_{\dff n}
\off,\qquad
p_{\fff j\dff -\dff 1}
\off =\off
l_{\dff 1}\dff \cdot\qff
l_{\trf 2}\dff \cdot\qff
\ldots\qff \cdot\qff
l_{\dff n}
\off.
\]

\vspace{-9pt}
Finally\halfff,\oss let\dss 
$M_{\trf 0}\dff,\off M_{\dff 1}\dff,\off \ldots\dff,\off M_{\dff n}$\dss
be\sss the vertical\sss segments of\trs the above form contained\dss in\dss
the rectangle\dss $Q_{\trf j}$\dss and\trs listed\dss from\dss left\dss to\sss right\halfff,\oss
and\trs let\dss 
$m_{\dff i}\off =\off h_{\qff M_{\dff i}}$.\oss
Then every\dss $m_{\dff i}$\dss is\dss a\sss short\dss path and
$m_{\qff 0}\dff,\off m_{\dff n}$ are constant\dss paths.\oss
By\dss the definition,\vspace{3pt}
\[
\quad
u^{\dff \sharp}\dff(\dff p_{\fff j}\trf)
\off =\off
u\trf(\dff k_{\dff 1}\trf)\dff \cdot\dff
u\trf(\dff k_{\qff 2}\trf)\dff \cdot\dff
\ldots\dff \cdot\dff
u\trf(\dff k_{\dff n}\trf)
\off,\qquad
u^{\dff \sharp}\dff(\dff p_{\fff j\dff -\dff 1}\trf)
\off =\off
u\trf(\dff l_{\dff 1}\trf)\dff \cdot\dff
u\trf(\dff l_{\dff 2}\trf)\dff \cdot\dff
\ldots\dff \cdot\dff
u\trf(\dff l_{\dff n}\trf)
\pff.
\]

\vspace{-9pt}
As we pointed out\sss at\sss the beginning of\qss Section\qss \ref{cochains},\pss
if\dss $m$\sss is\dss a constant\dss path,\oss
then\dss $u\trf(\dff m\trf)\off =\off 1$\trs
By\sss applying\dss this remark\dss to\dss
$m_{\qff 0}\dff,\off m_{\dff n}$\dss
we see\sss that\vspace{3pt}
\[
\quad
u^{\dff \sharp}\dff(\dff p_{\fff j}\trf)\hspace{1.15em}
\off =\off
u\trf(\dff m_{\qff 0}\trf)\dff \cdot\dff
u\trf(\dff k_{\dff 1}\trf)\dff \cdot\dff
u\trf(\dff k_{\qff 2}\trf)\dff \cdot\dff
\ldots\dff \cdot\dff
u\trf(\dff k_{\dff n}\trf)\quad\
\mbox{and}\quad\
\]

\vspace{-36pt}
\[
\quad
u^{\dff \sharp}\dff(\dff p_{\fff j\dff -\dff 1}\trf)
\off =\off
u\trf(\dff l_{\dff 1}\trf)\dff \cdot\dff
u\trf(\dff l_{\dff 2}\trf)\dff \cdot\dff
\ldots\dff \cdot\dff
u\trf(\dff l_{\dff n}\trf)
\dff \cdot\dff
u\trf(\dff m_{\dff n}\trf)
\pff.
\]

\vspace{-9pt}
Therefore,\oss it\dss is\dss sufficient\dss to 
prove\sss that\dss these\sss two products are equal.\oss

\vspace{-6pt}
\begin{figure}[h!]
\hspace*{10em}
\includegraphics[width=0.48\textwidth]{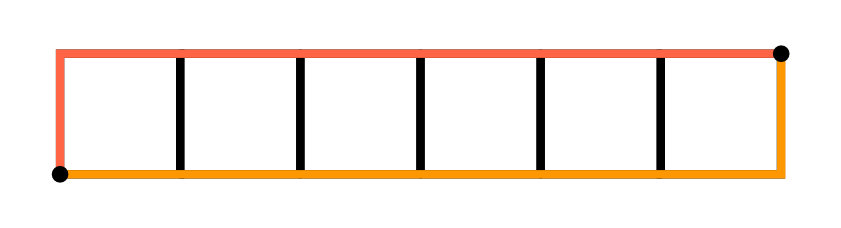}
\end{figure}\vspace{-6pt}

We will\dss prove\sss that\halfff,\oss moreover\halfff,\oss
all\dss products of\trs the form\vspace{3pt}
\[
\quad
g_{\dff i}
\off =\off
u\trf(\dff l_{\dff 1}\trf)\dff \cdot\dff
\ldots\dff \cdot\dff
u\trf(\dff l_{\dff i\dff -\dff 1}\trf)\dff \cdot\dff
u\trf(\dff m_{\dff i}\trf)\dff \cdot\dff
u\trf(\dff k_{\dff i}\trf)\dff \cdot\dff
\ldots\dff \cdot\dff
u\trf(\dff k_{\dff n}\trf)
\]

\vspace{-9pt}
with\dss
$i\off =\off 0\fff,\pff 1\fff,\pff \ldots\fff,\pff n$\dss
are equal.\oss

\vspace{-6pt}
\begin{figure}[h!]
\hspace*{10em}
\includegraphics[width=0.48\textwidth]{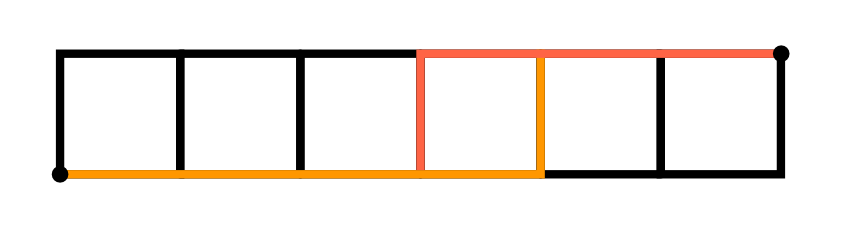}
\end{figure}\vspace*{-6pt}

\newpage
It\dss is\dss sufficient\dss to prove\sss that\dss
$g_{\dff i}\off =\off g_{\dff i\dff +\dff 1}$\dss
for\dss $i\qff \leq\qff n\qff -\qff 1$\nnsp.\oss
One\sss gets $g_{\dff i\dff +\dff 1}$\sss
by\sss replacing\dss in $g_{\dff i}$\sss
the\sss two consecutive factors
$u\trf(\dff m_{\dff i}\trf)\dff \cdot\dff
u\trf(\dff k_{\dff i}\trf)$\dss
by\dss the factors\dss 
$u\trf(\dff l_{\dff i}\trf)\dff \cdot\dff
u\trf(\dff m_{\dff i\dff +\dff 1}\trf)$\nnsp.\oss
Clearly\halfff,\oss the paths\dss
$m_{\dff i}\dff \cdot\dff k_{\dff i}$\dss
and\dss
$l_{\dff i}\dff \cdot\dff m_{\dff i\dff +\dff 1}$\dss
are short\sss and\dss homotopic\sss by\sss a\sss short\dss homotopy\halfff.\oss
It\dss follows\dss that\vspace{3pt}
\[
\quad
u\trf(\dff m_{\dff i}\trf)\dff \cdot\dff
u\trf(\dff k_{\dff i}\trf)
\off =\off
u\trf(\dff m_{\dff i}\dff \cdot\dff k_{\dff i}\trf)
\off =\off
u\trf(\dff l_{\dff i}\dff \cdot\dff m_{\dff i\dff +\dff 1}\trf)
\off =\off
u\trf(\dff l_{\dff i}\trf)\dff \cdot\dff
u\trf(\dff m_{\dff i\dff +\dff 1}\trf)
\]

\vspace{-9pt}
and\dss hence\dss
$g_{\dff i}\off =\off g_{\dff i\dff +\dff 1}$\nsp.\oss
As we saw,\oss this implies\sss that\trs
$u^{\dff \sharp}\dff(\dff p_{\fff j}\trf)
\off =\off 
u^{\dff \sharp}\dff(\dff p_{\fff j\dff -\dff 1}\trf)$\nnsp.\oss
In\dss turn,\oss
this implies\sss that\trs
$u^{\dff \sharp}\dff(\dff p\trf)
\off =\off 
u^{\dff \sharp}\dff(\dff q\trf)$\nnsp.\oss
This completes\sss the proof\dss of\trs the fact\dss that\dss $u^{\dff \sharp}$\dss
is\dss a cocycle.\oss

Since\sss
the cocycle\dss $u^{\dff \sharp}$\dss does not\sss depend
on\dss the partitions of\trs $I$\dss used\sss in\dss the construction,\oss 
the restriction of\dss $u^{\dff \sharp}$\dss
to\dss $P_{\dff \mathcal{U}}\dff(\trf X\trf)$\dss is\dss equal\dss to\sss
the original\dss short\dss cocycle $u$\dss
(for a short\dss path one can use\sss the partition of\trs $I$\dss
into one interval,\oss namely\dss $I$\nnsp).\oss
Conversely\halfff,\oss suppose\sss that\trs 
$v\qff \in\qff Z^{\dff 1}\dff(\trf X\fff,\pff Y\trf)$\dss
and\dss $u$\dss is\dss the restriction of\dss $v$\dss to\dss
$P_{\dff \mathcal{U}}\dff(\trf X\trf)$\nnsp.\oss
Since\dss $v$\dss is\dss a\sss cocycle,\oss
(\ref{u-sharp})\qss implies\sss that\dss
$u^{\dff \sharp}\off =\off v$\nnsp.\oss
It\dss follows\sss that\dss the restriction\dss map\sss 
$Z^{\dff 1}\dff(\trf X\fff,\pff Y\trf)
\off \ttoo\off
Z^{\dff 1}_{\dff \mathcal{U}}\dff(\trf X\fff,\pff Y\trf)$\sss
is\dss a\dss bijection.\oss

Clearly\halfff,\oss the restriction of\dss cocycles\sss to
$P_{\dff \mathcal{U}}\dff(\trf X\trf)$\sss
is\dss equivariant\dss with\dss respect\dss to\sss
the action of\sss $0$\dnsp-cochains\qss
(which are all\sss short\fff)\qss
on\sss $Z^{\dff 1}\dff(\trf X\fff,\pff Y\trf)$\sss
and\sss $Z^{\dff 1}_{\dff \mathcal{U}}\dff(\trf X\fff,\pff Y\trf)$\nnsp.\oss
Since\dss the above restriction map\dss is\dss a\sss bijection,\oss
the equivariance implies\sss that\dss the induced map\sss
$H^{\dff 1}\dff(\trf X\fff,\pff Y\trf)
\off \ttoo\off
H^{\dff 1}_{\dff \mathcal{U}}\dff(\trf X\fff,\pff Y\trf)$\sss
is\dss also\sss a\dss bijection.\oss
This completes\sss the proof\dss of\trs the\sss theorem.\oss  \eproof

\mysection{The\qss standard\qss Seifert--van Kampen\qss theorem}{svk-main}

\myuppar{The restriction maps.}
Suppose\sss that\sss $A$\sss is\dss
a path-connected subspace of\dss $X$\dss and\dss $b\qff \in\qff A$\nnsp.\oss
The inclusion map\sss $i\dff \colon\dff A\qff \ttoo\qff X$\dss
induces\sss the homomorphism\dss
$i_{\dff *}\dff \colon\dff
\pi_{\dff 1}\trf(\trf A\fff,\qff b\trf)
\qff \ttoo\qff 
\pi_{\dff 1}\trf(\trf X\fff,\qff b\trf)$\nnsp,\oss
which,\oss in\dss turn,\oss induces\sss the\qss
\emph{restriction\sss map}\vspace{3pt}\vspace{-0.375pt}
\[
\quad
r_{\qff A}\dff \colon\dff
\hom (\trf \pi_{\dff 1}\trf(\trf X\fff,\qff b\trf)\fff,\pff G\trf)
\off \ttoo\off
\hom (\trf \pi_{\dff 1}\trf(\trf A\fff,\qff b\trf)\fff,\pff G\trf)
\pff.
\]

\vspace{-9pt}\vspace{-0.375pt}
Similarly,\qss $i$ induces a map\dss
$i_{\dff *}\dff \colon\dff
P\dff(\trf A\trf)\qff \ttoo\qff P\dff(\trf X\trf)$
such\sss that\sss 
$i_{\dff *}\dff(\dff p\dff \cdot\dff q\trf)
\off =\off
i_{\dff *}\dff(\dff p\trf)
\dff \cdot\dff 
i_{\dff *}\dff(\dff q\trf)$\dss
when\dss $p\dff \cdot\dff q$\dss is\dss defined,\oss
and $i\fff,\qff i_{\dff *}$ induce\sss the\qss 
\emph{restriction\dss maps}\vspace{3pt}
\[
\quad
r_{\qff A}\dff \colon\dff
C^{\dff 0}\dff(\trf X\fff,\pff b\trf)
\off \ttoo\off
C^{\dff 0}\dff(\trf A\fff,\pff b\trf)
\dff,\quad
r_{\qff A}\dff \colon\dff
F^{\dff 1}\dff(\trf X\fff,\pff b\trf)
\off \ttoo\off
F^{\dff 1}\dff(\trf A\fff,\pff b\trf)
\pff
\]

\vspace{-9pt}
respectively,\oss
where\dss $F$\dss stands\dss for\sss $C$\nnsp,\dss $Z$\nnsp,\qss or\dss $H$\nnsp.\oss
Let\sss $\mathcal{U}$ be an open covering of\trs $X$\nnsp.\oss
Then\dss $\{\qff U\dff \cap\dff A\dff \mid\dff U\qff \in\qff \mathcal{U}\qff\}$\dss is\dss
an open covering of\trs $A$\nnsp,\oss
which we will\sss still\sss denote by\sss $\mathcal{U}$\nnsp.\oss
The inclusion\sss 
$i$\sss induces a map\dss
$i_{\dff *}\dff \colon\dff
P_{\dff \mathcal{U}}\dff(\trf A\trf)
\qff \ttoo\qff 
P_{\dff \mathcal{U}}\dff(\trf X\trf)$\nnsp,\oss
which,\oss in\dss turn,\oss induces\qss 
\emph{restriction\sss maps}\vspace{3pt}
\[
\quad
r_{\qff A}\dff \colon\dff
F^{\dff 1}_{\dff \mathcal{U}}\dff(\trf X\fff,\pff b\trf)
\off \ttoo\off
F^{\dff 1}_{\dff \mathcal{U}}\dff(\trf A\fff,\pff b\trf)
\pff,
\]

\vspace{-9pt}
where\dss $F\off =\off C\fff,\pff Z$\dss or\dss $H$\nnsp.\oss

\mypar{Theorem.}{svk-olum}
\emph{Suppose\sss that\trs $U\fff,\pff V\qff \subset\qff X$\dss
are\sss two path-connected open sets such\dss that\dss $U\qff \cap\qff V$\dss
is\dss path-connected and\dss $b\qff \in\qff U\qff \cap\qff V$\dnsp.\oss
Then\dss the square of\dss restriction maps}\vspace{-1.5pt}
\[
\quad
\begin{tikzcd}[column sep=large, row sep=normal]
& 
H^{\dff 1}\dff(\trf U\fff,\pff b\dff)
\arrow[rd, near start, "{\dis r_{\trf U\dff \cap\dff V}}"]
&\\
H^{\dff 1}\dff(\trf X\fff,\pff b\dff) 
\arrow[ru, "{\dis r_{\trf U}}"]
\arrow[rd, "{\dis r_{\trf V}}"'] 
&
&
H^{\dff 1}\dff(\trf U\dff \cap\dff V\fff,\pff b\dff)\\
&
H^{\dff 1}\dff(\trf V\fff,\pff b\dff)
\arrow[ru, near start, "{\dis r_{\trf U\dff \cap\dff V}}"']
&
\end{tikzcd}
\]

\vspace{-12pt}\vspace{-1.5pt}
\emph{is\dss commutative and cartesian.\oss}

\proof
Let\dss $\mathcal{U}\off =\off \{\trf U\fff,\pff V\trf\}$\nnsp.\oss
Let\dss us\dss replace\dss the cohomology\sss set\trs 
$H^{\dff 1}\dff(\trf X\fff,\pff b\dff)$\dss 
by\dss its\sss short\dss version\dss
$H^{\dff 1}_{\dff \mathcal{U}}\dff(\trf U\fff,\pff b\dff)$\nnsp.\oss
It\dss is\dss sufficient\dss to prove\sss that\dss 
the resulting square\vspace{-1.5pt}
\[
\quad
\begin{tikzcd}[column sep=large, row sep=normal]
& 
H^{\dff 1}\dff(\trf U\fff,\pff b\dff)
\arrow[rd, near start, "{\dis r_{\trf U\dff \cap\dff V}}"]
&\\
H^{\dff 1}_{\dff \mathcal{U}}\dff(\trf X\fff,\pff b\dff) 
\arrow[ru, "{\dis r_{\trf U}}"]
\arrow[rd, "{\dis r_{\trf V}}"'] 
&
&
H^{\dff 1}\dff(\trf U\dff \cap\dff V\fff,\pff b\dff)\pff,\\
&
H^{\dff 1}\dff(\trf V\fff,\pff b\dff)
\arrow[ru, near start, "{\dis r_{\trf U\dff \cap\dff V}}"']
&
\end{tikzcd}
\]

\vspace{-12pt}\vspace{-1.5pt}
where all\dss maps are\sss the restriction maps,\oss
is\dss commutative and cartesian.\oss
The commutativity\dss is\dss obvious
both\dss for\dss the original\sss square and\dss its short\dss version.\oss
To simplify\dss notations,\oss
we will\sss omit\dss the base point\dss $b$\sss in\dss the rest\sss of\trs the proof\halfff.\oss
The\sss last\sss square\dss is\dss cartesian\dss
if\trs the map\vspace{0.75pt}
\[
\quad
(\dff r_{\trf U}\fff,\pff r_{\trf V}\trf)\dff \colon\dff
H^{\dff 1}_{\dff \mathcal{U}}\dff(\trf X\trf)
\off \ttoo\off
H^{\dff 1}\dff(\trf U\trf) 
\qff \times\qff
H^{\dff 1}\dff(\trf V\trf) 
\]

\vspace{-12pt}\vspace{0.75pt}
induces a bijection\sss from\sss 
$H^{\dff 1}_{\dff \mathcal{U}}\dff(\trf X\trf)$\dss
to\sss the fibered\dss product\sss 
$H^{\dff 1}\dff(\trf U\trf) 
\qff \times_{\pff H^{\dff 1}\dff(\trf U\qff \cap\qff V\trf)}\qff
H^{\dff 1}\dff(\trf V\trf)$\nnsp.\oss

\emph{Surjectivity\halfff.}\oss 
Suppose\sss that\trs
$u\qff \in\qff Z^{\dff 1}\dff(\trf U\trf)$\dss
and\dss
$v\qff \in\qff Z^{\dff 1}\dff(\trf V\trf)$\dss
are such\dss that
\[
\quad
r_{\trf U\dff \cap\dff V}\dff(\dff u\trf)\quad\
\mbox{and}\quad\
r_{\trf U\dff \cap\dff V}\dff(\dff v\trf)
\]

\vspace{-12pt}
belong\dss to\sss the same cohomology\sss class.\oss
Then\dss there exists a $0$\dnsp-cochain\dss
$c\qff \in\qff C^{\dff 0}\dff(\trf U\qff \cap\qff V\trf)$\dss
such\dss that\trs 
$c\dff \bullet\dff v\dff(\dff p\trf)
\off =\off
u\dff(\dff p\trf)$\dss
for every\dss path\dss $p\qff \in\qff P\trf(\trf U\qff \cap\qff V\trf)$\nnsp.\oss
Let\dss us extend\dss the cochain\dss $c$\dss to a cochain\dss
$c^{\dff \sim}\qff \in\qff C^{\dff 0}\dff(\trf V\trf)$\nnsp.\oss
The cochains\sss $u$\sss and\dss 
$c^{\dff \sim}\dff \bullet\dff v$\dss
agree on\dss
$P\trf(\trf U\qff \cap\qff V\trf)
\off =\off
P\trf(\trf U\trf)
\qff \cap\qff
P\trf(\trf V\trf)$\dss
and\dss hence define a $1$\dnsp-cochain\dss
$w\qff \in\qff C^{\dff 1}_{\dff \mathcal{U}}\dff(\trf X\trf)$\nnsp.\oss
The cochain\sss $w$\sss is\dss a\sss short\sss cocycle\sss because\sss the conditions for beings
a short\sss cocycle are imposed only\sss in\sss $U$\sss and\dss in\sss $V$\dnsp.\oss
Clearly\halfff,\oss\vspace{0.75pt}
\[
\quad
r_{\trf U}\dff(\dff w\trf)
\off =\off
u\quad\
\mbox{and}\quad\
r_{\trf V}\dff(\dff w\trf)
\off =\off
c^{\dff \sim}\dff \bullet\dff v
\pff.
\]

\vspace{-12pt}\vspace{0.75pt}
Since\dss $c^{\dff \sim}\dff \bullet\dff v$\dss belongs\sss to\sss the same cohomology\sss
class as\sss $v$\nnsp,\oss the surjectivity\dss follows.\oss

\emph{Injectivity\halfff.}\oss 
Suppose\sss that\trs
$w\halfff,\pff z\pff \in\pff Z^{\dff 1}_{\dff \mathcal{U}}\dff(\trf X\trf)$\dss
are such\dss that\trs
$r_{\trf U}\dff(\dff w\trf)$\dss
and\dss
$r_{\trf U}\dff(\dff z\trf)$\dss
belong\sss to\sss the same cohomology\sss class\sss in\dss
$H^{\dff 1}\dff(\trf U\trf)$\nnsp,\oss
and\dss 
$r_{\trf V}\dff(\dff w\trf)$\dss
and\dss
$r_{\trf V}\dff(\dff z\trf)$\dss
belong\sss to\sss the same cohomology\sss class\sss in\dss
$H^{\dff 1}\dff(\trf V\trf)$\nnsp.\oss
Then\dss there exist\sss $0$\dnsp-cochains\dss
$a\qff \in\qff C^{\dff 0}\dff(\trf U\trf)$\dss
and\dss
$c\qff \in\qff C^{\dff 0}\dff(\trf V\trf)$\dss
such\dss that\vspace{3pt}
\[
\quad
w\trf(\dff p\trf)
\off =\off
a\trf(\trf p\dff(\dff 0\dff)\trf)\dff \cdot\dff
z\trf(\dff p\trf)\dff \cdot\dff
a\trf(\trf p\dff(\dff 1\dff)\trf)^{\dff -\dff 1}
\quad\
\mbox{and}
\]

\vspace{-36pt}
\[
\quad
w\trf(\dff p\trf)
\off =\off
c\trf(\trf p\dff(\dff 0\dff)\trf)\dff \cdot\dff
z\trf(\dff p\trf)\dff \cdot\dff
c\trf(\trf p\dff(\dff 1\dff)\trf)^{\dff -\dff 1}
\pff
\]

\vspace{-9pt}
if\dss $p$\dss belongs\sss to\dss $P\trf(\trf U\trf)$\dss
and\dss $P\trf(\trf V\trf)$\dss respectively\halfff.\oss
We claim\dss that\dss 
$a\dff(\dff x\trf)\off =\off c\trf(\dff x\trf)$\dss
for every\dss $x\qff \in\qff U\qff \cap\qff V$\nnsp.\oss
Indeed,\oss since\dss $U\qff \cap\qff V$\dss
is\dss path-connected,\oss
there exists\dss 
$p\qff \in\qff P\trf(\trf U\qff \cap\qff V\trf)$\dss
such\dss that\dss
$p\dff(\dff 0\dff)\off =\off b$\dss
and\dss
$p\dff(\dff 1\dff)\off =\off x$\nnsp.\oss
Then\dss
$a\trf(\trf p\dff(\dff 0\dff)\trf)
\off =\off
c\trf(\trf p\dff(\dff 0\dff)\trf)
\off =\off
1$\nnsp.\oss
It\dss follows\dss that\vspace{2.875pt}
\[
\quad
z\trf(\dff p\trf)\dff \cdot\dff
a\trf(\trf x\trf)^{\dff -\dff 1}
\off =\off
z\trf(\dff p\trf)\dff \cdot\dff
c\trf(\trf x\trf)^{\dff -\dff 1}
\]

\vspace{-9.125pt}
and\dss hence\dss
$a\trf(\trf x\trf)
\off =\off
c\trf(\trf x\trf)$\dss
for every\dss $x\qff \in\qff U\qff \cap\qff V$\dnsp.\oss
In other\dss terms,\pss $a$\sss and\sss $c$\sss agree on\dss the 
intersection of\dss their domains and\dss hence define a 
$0$\dnsp-cochain\dss $d\qff \in\qff C^{\dff 0}\dff(\trf X\fff,\pff b\dff)$\dss
such\dss that\vspace{2.875pt}
\[
\quad
w\trf(\dff p\trf)
\off =\off
d\trf(\trf p\dff(\dff 0\dff)\trf)\dff \cdot\dff
z\trf(\dff p\trf)\dff \cdot\dff
d\trf(\trf p\dff(\dff 1\dff)\trf)^{\dff -\dff 1}
\]

\vspace{-9.125pt}
for every\sss short\dss path\sss $p$\nnsp.\oss
Therefore\sss the cohomology\sss classes of\trs $w$\sss and\dss $z$\dss
in\dss $H^{\dff 1}_{\dff \mathcal{U}}\dff(\trf X\trf)$\dss
are equal.\oss
The injectivity\dss follows.\oss
This completes\sss the proof\dss of\trs the\sss theorem.\oss  \eproof

\mypar{Theorem.}{svk-universal}
\emph{Under\dss the assumptions of\qss Theorem\qss \ref{svk-olum}\qss
the fundamental\dss group\dss 
$\pi_{\dff 1}\trf(\trf X\fff,\qff b\trf)$\dss
has\sss the following\sss universal\dss property\halfff.\oss
Let\dss us\dss consider\dss diagrams of\trs the form}\vspace{6pt}\vspace{2pt}
\[
\quad
\begin{tikzcd}[column sep=normal, row sep=hugeplus]
& 
\pi_{\dff 1}\dff(\trf U\fff,\pff b\dff)
\arrow[d]
\arrow[rrd, "{\dis h_{\trf U}}"]
&
&\\
\pi_{\dff 1}\dff(\trf U\dff \cap\dff V\fff,\pff b\dff) 
\arrow[ru]
\arrow[rd] 
&
\pi_{\dff 1}\dff(\trf X\fff,\pff b\dff)
\arrow[rr, dashed]
&
&
G\pff,\phantom{GGG}\\
&
\pi_{\dff 1}\dff(\trf V\fff,\pff b\dff)
\arrow[u]
\arrow[rru, "{\dis h_{\trf V}}"']
&
&
\end{tikzcd}
\]

\vspace{-6pt}\vspace{2pt}
\emph{where all\dss unmarked solid arrows are\sss
homomorphisms induced\dss by\dss inclusions,\pss
$G$\dss is\dss a\sss group,\oss
and\qss $h_{\trf U}\dff,\pff h_{\trf V}$\dss
are homomorphisms such\dss that\dss the outer square\dss is\dss commutative.\oss
For every\dss such diagram\dss there exists unique dashed arrow
making\dss the\sss two right\dss triangles commutative.\oss}

\proof
To\sss simplify\dss notations,\pss
we will\sss again omit\dss the base point\dss $b$\nnsp.\oss
In\dss the\sss language of\trs the sets\dss
$\hom (\trf \pi_{\dff 1}\trf(\trf W\trf)\fff,\pff G\trf)$\nnsp,\oss
where\dss $W\off =\off X\fff,\pff U\fff,\pff V$ or\dss $U\qff \cap\qff V$\dnsp,\oss
the\sss theorem\dss claims\sss that\dss the square\vspace{6pt}\vspace{2pt}
\[
\quad
\hspace*{-2.5em}
\begin{tikzcd}[column sep=small, row sep=huge]
& 
\hom (\trf \pi_{\dff 1}\trf(\trf U\trf)\fff,\pff G\trf)
\arrow[rd, "{\dis r_{\trf U\dff \cap\dff V}}"]
&\\
\phantom{GGg}\hom (\trf \pi_{\dff 1}\trf(\trf X\trf)\fff,\pff G\trf)\phantom{GGg}
\arrow[ru, "{\dis r_{\trf U}}"]
\arrow[rd, "{\dis r_{\trf V}}"'] 
&
&
\hom (\trf \pi_{\dff 1}\trf(\trf U\qff \cap\qff V\trf)\fff,\pff G\trf)\phantom{GGGG}\\
&
\hom (\trf \pi_{\dff 1}\trf(\trf V\trf)\fff,\pff G\trf)
\arrow[ru, "{\dis r_{\trf U\dff \cap\dff V}}"']
&
\end{tikzcd}
\]

\vspace{-6pt}\vspace{2pt}
is\dss cartesian.\oss
This follows\sss from\trs
Theorem\qss \ref{svk-olum},\oss
the identification of\trs the sets\dss
$\hom (\trf \pi_{\dff 1}\trf(\trf Z\trf)\fff,\pff G\trf)$\dss
with\dss the sets\dss
$H^{\dff 1}\dff(\trf Z\fff,\pff b\dff)$\dss
as in\dss Section\qss \ref{cochains},\oss
and\dss the fact\dss that\dss this identification
obviously agrees with\dss the restriction maps.\oss  \eproof

\myuppar{Theorem\qss \ref{svk-universal},\oss free products,\oss and\sss relations.}
For\sss $W\off =\off U$\sss or\dss $V$\sss let\sss
$i\dff(\trf W\trf)\dff \colon\dff
U\qff \cap\qff V\qff \ttoo\qff W$\sss
be\sss the inclusion map.\oss
The universal\dss property\dss of\pss Theorem\qss \ref{svk-universal}\qss
means\sss that\sss $\pi_{\dff 1}\dff(\trf X\fff,\pff b\dff)$\sss
is\dss iso\-mor\-phic\sss to\sss the free product\sss
$\pi_{\dff 1}\dff(\trf U\fff,\pff b\dff)\dff *\dff \pi_{\dff 1}\dff(\trf V\fff,\pff b\dff)$\dss
with\sss the relation\dss 
$i\dff(\trf U\trf)_{\dff *}\dff(\trf \gamma\trf)
\off =\off
i\dff(\trf V\trf)_{\dff *}\dff(\trf \gamma\trf)$\sss 
imposed\dss for every\sss
$\gamma\qff \in\qff \pi_{\dff 1}\dff(\trf U\qff \cap\qff V\fff,\pff b\dff)$\nnsp.\oss

\mysection{Unions\qss of\qss several\qss subsets}{svk-families}

\myuppar{Open coverings.}
Theorems\qss \ref{svk-olum}\qss and\qss \ref{svk-universal}\qss
can\sss be\sss generalized\dss to\sss the following\dss situation.\oss
Let\sss $\mathcal{U}$\dss be an open covering of\sss $X$\nnsp.\oss
Suppose\sss that\trs $b\qff \in\qff U$\dss for every\dss 
$U\qff \in\qff \mathcal{U}$\dnsp,\oss
that\sss every\sss set\sss $U\qff \in\qff \mathcal{U}$\dss
is\dss path-connected,\oss
the intersection\dss $U\qff \cap\qff V$\dss
is\dss path-connected for every\dss two sets\sss 
$U\fff,\pff V\qff \in\qff \mathcal{U}$\nnsp,\oss
and\dss the intersection\sss $U\qff \cap\qff V\qff \cap\qff W$\sss
is\dss path-connected for every\dss three sets\dss
$U\fff,\pff V\fff,\pff W\qff \in\qff \mathcal{U}$\nnsp.\oss

\mypar{Theorem.}{olum-family}
\emph{Suppose\sss that\dss the above assumptions hold.\oss 
Given a\sss family}\qss \vspace{3pt}
\[
\quad
\left\{\qff
h_{\trf U}\qff \in\qff H^{\dff 1}\dff(\trf U\fff,\pff b\dff)
\qff\right\}_{\qff U\qff \in\qff \mathcal{U}}
\qff,
\]

\vspace{-9pt}
\emph{there exists\sss
$h\qff \in\qff H^{\dff 1}\dff(\trf X\fff,\pff b\dff)$\sss 
such\dss that\sss
$r_{\trf U}\trf(\dff h\trf)
\off =\off 
h_{\trf U}$\sss
for every\sss $U\qff \in\qff \mathcal{U}$\sss
if\trs and\dss only\trs if}\vspace{3pt}
\[
\quad
r_{\trf U\dff \cap\dff V}\trf(\dff h_{\trf U}\trf)
\off =\off 
r_{\trf U\dff \cap\dff V}\trf(\dff h_{\trf V}\trf)
\]

\vspace{-9pt}
\emph{for every\sss pair\sss $U\fff,\pff V\qff \in\qff \mathcal{U}$\dnsp.\oss
If\dss such a cohomology\sss class $h$ exists,\oss it\dss is\dss unique.\oss}

\proof
As\sss in\dss the proof\dss of\qss Theorem\qss \ref{svk-olum},\oss
we can\dss replace\dss 
$H^{\dff 1}\dff(\trf X\fff,\pff b\dff)$\dss 
by\dss 
$H^{\dff 1}_{\dff \mathcal{U}}\dff(\trf U\fff,\pff b\dff)$\nnsp.\oss
We will\sss omit\dss the base point\dss $b$\sss in\dss the rest\sss of\dss the proof\halfff.\oss
Obviously\halfff,\oss the stated condition\dss is\dss necessary\halfff.

Let\dss us\dss prove\sss first\dss that\dss if\sss 
$h$\sss exists,\oss then\dss it\dss is\dss unique.\oss
Suppose\sss that\sss
$w\halfff,\pff z\pff \in\pff Z^{\dff 1}_{\dff \mathcal{U}}\dff(\trf X\trf)$\sss
are such\dss that\sss
$r_{\trf U}\dff(\dff w\trf)$\sss
and\sss
$r_{\trf U}\dff(\dff z\trf)$\sss
belong\sss to\sss the same cohomology\sss class\sss in\sss
$H^{\dff 1}\dff(\trf U\trf)$\sss
for every\sss $U\qff \in\qff \mathcal{U}$\dnsp.\oss
Then\dss for every\sss $U\qff \in\qff \mathcal{U}$\sss
there exists a $0$\dnsp-cochain\sss
$c_{\qff U}\qff \in\qff C^{\dff 0}\dff(\trf U\trf)$\sss
such\dss that\vspace{3pt}
\[
\quad
w\trf(\dff p\trf)
\off =\off
c_{\qff U}\trf(\trf p\dff(\dff 0\dff)\trf)\dff \cdot\dff
z\trf(\dff p\trf)\dff \cdot\dff
c_{\qff U}\trf(\trf p\dff(\dff 1\dff)\trf)^{\dff -\dff 1}
\pff
\]

\vspace{-12pt}\vspace{3pt}
if\sss $p$\dss belongs\sss to\sss $P\trf(\trf U\trf)$\nnsp.\oss
By\sss using exactly\dss the same argument\sss as in\sss the proof\dss
of\dss Theorem\qss \ref{svk-olum},\oss
we see\sss that\sss
$c_{\qff U}\trf(\dff x\trf)\off =\off c_{\qff V}\trf(\dff x\trf)$\sss
for every\dss pair\sss $U\fff,\pff V\qff \in\qff \mathcal{U}$\sss
and\sss
$x\qff \in\qff U\qff \cap\qff V$\nnsp.\oss
It\dss follows\sss that\sss there exists\sss
$c\qff \in\qff C^{\dff 0}\dff(\trf X\fff,\pff b\dff)$\sss
such\dss that\sss
$c\trf(\dff x\trf)\off =\off c_{\qff U}\trf(\trf x\trf)$\sss
when\sss
$x\qff \in\qff U\qff \in\qff \mathcal{U}$\dnsp.\oss
Clearly\halfff,\oss\vspace{3pt}
\[
\quad
w\trf(\dff p\trf)
\off =\off
c\trf(\trf p\dff(\dff 0\dff)\trf)\dff \cdot\dff
z\trf(\dff p\trf)\dff \cdot\dff
c\trf(\trf p\dff(\dff 1\dff)\trf)^{\dff -\dff 1}
\]

\vspace{-12pt}\vspace{3pt}
for every\sss $p\qff \in\qff P_{\dff \mathcal{U}}\dff(\trf X\trf)$\nnsp.\oss
This proves\sss the uniqueness part\sss of\trs the\sss theorem.\oss

Suppose\sss now\dss that\dss for every\sss $U\qff \in\qff \mathcal{U}$\sss
a\sss $1$\dnsp-cocycle\sss
$u_{\trf U}\qff \in\qff Z^{\dff 1}\dff(\trf U\trf)$\sss
is\dss given and\dss that\sss 
$r_{\trf U\dff \cap\dff V}\dff(\dff u_{\trf U}\trf)$\sss 
and\sss
$r_{\trf U\dff \cap\dff V}\dff(\dff u_{\trf V}\trf)$\sss
belong\dss to\sss the same cohomology\sss class\sss
for every\sss $U\fff,\pff V\qff \in\qff \mathcal{U}$\dnsp.\oss
It\dss is\dss sufficient\dss to show\dss that\dss
there exist $1$\dnsp-cocycles\dss
$z_{\qff U}\qff \in\qff Z^{\dff 1}\dff(\trf U\trf)$\sss
such\sss that\sss $u_{\trf U}$\sss and\sss $z_{\qff U}$\sss
belong\sss to\sss the same cohomology\sss class 
for every\sss $U\qff \in\qff \mathcal{U}$\sss and\sss
$r_{\trf U\dff \cap\dff V}\dff(\dff z_{\qff U}\trf)
\off =\off
r_{\trf U\dff \cap\dff V}\dff(\dff z_{\qff V}\trf)$\sss
for every\sss $U\fff,\pff V\qff \in\qff \mathcal{U}$\dnsp.\oss
Indeed,\oss then\sss the cocycles\sss $z_{\qff U}$\sss
agree on\sss the intersections of\dss their domains and\dss hence define a cocycle\sss
$z\qff \in\qff Z_{\dff \mathcal{U}}^{\dff 1}\dff(\trf X\trf)$\nnsp.\oss
Let\sss $h$\sss be\sss the cohomology\sss class of\sss $z$\nnsp.\oss
Then\sss $r_{\trf U}\dff(\dff h\trf)$\sss is\dss equal\sss to\sss the cohomology\sss
class of\sss $u_{\trf U}$\sss for every\sss $U\qff \in\qff \mathcal{U}$\dnsp.\oss

Suppose\sss first\dss that\dss  
$\mathcal{U}$\dss 
is\dss finite,\oss say\sss 
$\mathcal{U}
\off =\off
\{
\qff
U_{\dff 1}\dff,\off
U_{\dff 2}\dff,\off
\ldots\dff,\off
U_{\dff m}
\qff
\}$\sss
for some\sss $m$\nnsp.\oss
Let\dss $u_{\dff i}\off =\off u_{\trf U_{\dff i}}$\nsp.\oss
Let\sss us assume\dss that\dss
there are cocycles\sss
$z_{\trf i}\qff \in\qff Z^{\dff 1}\dff(\trf U_{\dff i}\trf)$\dss
with\sss $i\qff \leq\qff m\qff -\qff 1$\sss
such\dss that\sss $z_{\trf i}$\sss belongs\sss to\sss the same cohomology\sss
class as\sss $u_{\dff i}$\sss for every\sss $i$\sss
and\dss the cocycles\sss $z_{\trf i}$\sss
agree on\dss the intersections of\dss their domains.\oss 
Then\dss for each\sss 
$i\qff \leq\qff m\qff -\qff 1$\sss 
there\dss exists\sss 
$c_{\trf i}
\qff \in\qff 
C^{\dff 0}\dff(\trf U_{\fff m}
\qff \cap\qff 
U_{\dff i}\trf)$\sss
such\dss that\trs\vspace{3pt}
\begin{equation}
\label{uz-cochain}
\quad
c_{\trf i}\trf \bullet\dff u_{\dff m}\trf(\dff p\trf)
\off =\off
z_{\trf i}\trf(\dff p\trf)
\end{equation}

\vspace{-9pt}
for every\dss 
$p\qff \in\qff P\trf(\trf U_{\fff m}\qff \cap\qff U_{\dff i}\trf)$\nnsp.\oss
Let\dss 
$x
\qff \in\qff
(\trf U_{\fff m}\qff \cap\qff U_{\dff i}\trf)
\qff \cap\qff
(\trf U_{\fff m}\qff \cap\qff U_{\dff j}\trf)
\off =\off
U_{\fff m}\qff \cap\qff U_{\dff i}\qff \cap\qff U_{\dff j}$\nsp.\oss
Since\dss $U_{\fff m}\qff \cap\qff U_{\dff i}\qff \cap\qff U_{\dff j}$\dss
is\dss assumed\dss to be path-connected,\oss
there exists\dss
$p\qff \in\qff P\trf(\trf U_{\fff m}\qff \cap\qff U_{\dff i}\qff \cap\qff U_{\dff j}\trf)$\dss
such\dss that\trs
$p\dff(\dff 0\dff)\off =\off b$\dss
and\dss
$p\dff(\dff 1\dff)\off =\off x$\nnsp.\oss
Then\dss
$z_{\trf i}\trf(\dff p\trf)\off =\off z_{\dff j}\trf(\dff p\trf)$\dss
and\pss (\ref{uz-cochain})\pss implies\sss that\vspace{3pt}
\[
\quad
z_{\trf i}\dff(\dff p\trf)
\off =\off
c_{\dff i}\dff(\trf p\dff(\dff 0\dff)\trf)\dff \cdot\dff
u_{\dff m}\trf(\dff p\trf)\dff \cdot\dff
c_{\dff i}\dff(\trf p\dff(\dff 1\dff)\trf)^{\dff -\dff 1}
\quad
\mbox{and}
\]

\vspace{-36pt}
\[
\quad
z_{\dff j}\trf(\dff p\trf)
\off =\off
c_{\dff j}\dff(\trf p\dff(\dff 0\dff)\trf)\dff \cdot\dff
u_{\dff m}\trf(\dff p\trf)\dff \cdot\dff
c_{\dff j}\dff(\trf p\dff(\dff 1\dff)\trf)^{\dff -\dff 1}
\pff.
\]

\vspace{-9pt}
Since\dss
$c_{\dff i}\dff(\trf p\dff(\dff 0\dff)\trf)
\off =\off 
c_{\dff i}\dff(\trf b\trf)
\off =\off 
1$\dss
and\dss
$c_{\dff j}\dff(\trf p\dff(\dff 0\dff)\trf)
\off =\off 
c_{\dff j}\dff(\trf b\trf)
\off =\off 
1$\nnsp,\oss
it\dss follows\dss that\vspace{3pt}
\[
\quad
c_{\dff i}\dff(\trf x\trf)
\off =\off 
c_{\dff i}\dff(\trf p\dff(\dff 0\dff)\trf)
\off =\off 
c_{\dff j}\dff(\trf p\dff(\dff 0\dff)\trf)
\off =\off
c_{\dff j}\dff(\trf x\trf)
\pff.
\]

\vspace{-9pt}
Therefore $0$\dnsp-cochains $c_{\dff i}$ agree on\dss the intersections 
and\sss define a $0$\dnsp-cochain\vspace{0pt}
\[
\quad
c\off \in\off
C^{\dff 0}\dff\left(\qff
U_{\fff m}\off \cap\qff \bigcup_{i\qff \leq\qff m\qff -\qff 1}\qff U_{\dff i}
\qff\right)
\pff.
\]

\vspace{-12pt}
Let\dss us\sss extend\sss $c$\sss to a $0$\dnsp-cochain\dss
$c^{\dff \sim}\qff \in\qff C^{\dff 0}\dff(\trf U_{\fff m}\trf)$\dss
and\dss set\dss
$z_{\trf m}
\off =\off
c^{\dff \sim}\dff \cdot\dff u_{\dff m}$\nsp.\oss
Clearly\halfff,\pss
$z_{\trf m}$\dss belongs\sss to\sss the same cohomology\sss class\sss as\sss
$u_{\dff m}$\sss
and\trs $z_{\trf m}\dff,\pff z_{\trf i}$\dss
agree on\dss the intersections of\trs their\sss domains for every\dss
$i\qff \leq\qff m\qff -\qff 1$\nnsp.\oss
An\sss induction\sss by $m$ completes\sss the proof\dss for finite $\mathcal{U}$\dnsp.\oss

This proves\sss the\sss theorem\sss for\sss finite families\dss $\mathcal{U}$\nnsp.\oss
Suppose now\dss that\dss $\mathcal{U}$\dss is\dss a countable family,\oss
say\sss
$\mathcal{U}
\off =\off
\{\qff
U_{\dff 1}\dff,\off
U_{\dff 2}\dff,\off
\ldots\dff,\off
U_{\dff i}\dff,\off
\ldots
\off\}$\nnsp.\oss
The construction of\trs the cocycle\sss $z_{\trf m}$\dss
in\dss the above proof\trs keeps\sss the already\sss constructed
cocycles\dss $z_{\trf i}$\dss with\dss $i\qff \leq\qff m\qff -\qff 1$\dss
intact\halfff.\oss
Therefore\sss this construction can\sss be continued\sss indefinitely\sss
and\dss leads\sss to a sequence of\dss cocycles\dss
$z_{\trf 1}\dff,\off
z_{\trf 2}\dff,\off
\ldots\dff,\off
z_{\trf i}\dff,\off
\ldots$\dss
such\dss that\dss $z_{\trf i}$\dss and\dss $u_{\dff i}$\dss
belong\dss to\sss the same cohomology\sss class for every\sss $i$\sss
and\dss the cocycles\dss $z_{\trf i}$\dss agree on\dss intersections
of\trs their domains.\oss
This proves\sss the\sss theorem\dss for countable families.\oss
In\dss fact\halfff,\oss the same argument\dss works for an arbitrary\dss family\dss 
$\mathcal{U}
\off =\off
\{
\qff
U_{\dff i}
\qff\}_{\qff i\qff \in\qff J}$\nsp.\oss
One only\dss needs\sss to\sss well-order\sss $J$\sss
and apply\dss the\sss transfinite induction,\oss
or use\qss Zorn\qss lemma.\oss  \eproof

\mypar{Theorem.}{svk-family-universal}
\emph{Under\dss the same assumptions,\oss 
$\pi_{\dff 1}\trf(\trf X\fff,\qff b\trf)$\dss
has\sss the following\sss universal\dss pro\-per\-ty\halfff.\oss
Suppose\sss that\dss $G$\dss is\dss a\sss group\dss and\qss 
$h_{\trf U}\dff \colon\dff
\pi_{\dff 1}\dff(\trf U\fff,\pff b\dff)
\qff \ttoo\qff
G$\dss 
is\dss a\sss homomorphism\dss for every\trs
$U\qff \in\qff \mathcal{U}$\dnsp.\oss 
If\pss
for\dss every\qss $U\fff,\pff V\qff \in\qff \mathcal{U}$\trs 
the square}\vspace{3pt}
\[
\quad
\hspace*{4em}
\begin{tikzcd}[column sep=large, row sep=large]
& 
\pi_{\dff 1}\dff(\trf U\fff,\pff b\dff)
\arrow[rd, "{\dis h_{\trf U}}"]
&\\
\pi_{\dff 1}\dff(\trf U\dff \cap\dff V\fff,\pff b\dff) 
\arrow[rd]
\arrow[ru]
&
&
\phantom{\pi_{\dff 1}\dff}G\phantom{(\trf\dff \cap\dff V\fff,\pff b\dff)}\\
&
\pi_{\dff 1}\dff(\trf V\fff,\pff b\dff)
\arrow[ru, "{\dis h_{\trf V}}"']
&
\end{tikzcd}
\]

\vspace{-12pt}\vspace{3pt}
\emph{is\dss commutative,\oss
then\dss there exists a unique homomorphism\dss
$h\dff \colon\dff
\pi_{\dff 1}\dff(\trf X\fff,\pff b\dff)\qff \ttoo\qff G$\dss
such\dss that}\vspace{3pt}
\[
\quad
\hspace*{4em}
\hspace*{2.2em}
\begin{tikzcd}[column sep=large, row sep=huge]
\pi_{\dff 1}\dff(\trf X\fff,\pff b\dff)
\arrow[rr, "{\dis h}"]
\arrow[rd]
&
&
G\\
&
\pi_{\dff 1}\dff(\trf U\fff,\pff b\dff)
\arrow[ru, "{\dis h_{\trf U}}"']
&
\end{tikzcd}
\]

\vspace{-12pt}\vspace{3pt}
\emph{is\dss a\sss commutative\sss triangle\qss 
for\dss every\qss $U\qff \in\qff \mathcal{U}$\dnsp.\oss}

\proof
{\dff}It\trs is\dss completely\sss similar\dss to\sss the proof\dss of\qss
Theorem\qss \ref{svk-universal}.\oss  \eproof

\myuppar{Theorem\qss \ref{svk-family-universal},\oss free products,\oss and\sss relations.}
Given\sss $U\fff,\pff V\qff \in\qff \mathcal{U}$\dnsp,\oss 
let\trs \vspace{1.5pt}
\[
\quad
i\dff(\trf U\dff \cap\dff V\fff,\pff U\trf)\dff \colon\dff
U\qff \cap\qff V\qff \ttoo\qff U
\]

\vspace{-12pt}\vspace{1.5pt}
be\sss the inclusion map.\oss
The universal\dss property\dss of\pss Theorem\qss \ref{svk-family-universal}\qss
means\sss that\trs $\pi_{\dff 1}\dff(\trf X\fff,\pff b\dff)$\dss
is\dss iso\-mor\-phic\sss to\sss the free product\sss of\dss groups\sss
$\pi_{\dff 1}\dff(\trf U\fff,\pff b\dff)$\nnsp,\dss
$U\qff \in\qff \mathcal{U}$\dnsp,\oss
with\dss the\dss relations\dss\vspace{1.5pt}
\[
\quad
i\dff(\trf U\dff \cap\dff V\fff,\pff U\trf)_{\dff *}\dff(\trf \gamma\trf)
\off =\off
i\dff(\trf U\dff \cap\dff V\fff,\pff V\trf)_{\dff *}\dff(\trf \gamma\trf) 
\]

\vspace{-12pt}\vspace{1.5pt}
imposed\dss for every\sss
$U\fff,\pff V\qff \in\qff \mathcal{U}$\sss
and\sss
$\gamma
\qff \in\qff 
\pi_{\dff 1}\dff(\trf U\qff \cap\qff V\fff,\pff b\dff)$\nnsp.\oss\vspace{-1pt}

\myuppar{Remarks.}
The assumption of\dss the path-connectedness of\dss
triple intersections\sss $U\qff \cap\qff V\qff \cap\qff W$\sss
was used only\sss in\sss the paragraph\sss following\sss
the formula\qss (\ref{uz-cochain})\qss in\sss the proof\dss 
of\trs Theorem\qss \ref{olum-family}.\oss

\mysection{van Kampen\qss theorems}{relative-cohomology}

\myuppar{The case of\dss two point\sss subset $Y$\dnsp.}
The main\sss goal\sss of\dss this section\dss is\dss to
extend\sss the results of\trs Section\qss \ref{svk-main}\qss
to\sss the situation when\sss the intersection $U\dff \cap\dff V$\sss
is\dss not\sss necessarily\sss path-connected.\oss
The main ideas are present\sss already in\sss the case when\sss
this intersection consists of\dss two path-connected components,\oss
and we discuss\sss this case first.\oss
This requires some preliminary discussion of\dss
$H^{\dff 1}\dff(\trf X\dff,\pff Y\trf)$
with\sss $Y$ consisting of\dss two points\qss
(and,\oss in particular\halfff,\oss discrete).\oss

Recall\sss that\sss $b\qff \in\qff Y$ and\dss
let\sss us assume\sss that\sss $Y\off =\off \{\dff a\fff,\qff b\trf\}$
for some $a\qff \neq\qff b$\nnsp.\oss
Let\sss $G_{\dff a}\qff \subset\qff C^{\dff 0}\dff(\trf X\dff,\pff b\trf)$
be\sss the subgroup of $0$\dnsp-cochains equal\sss to $1\qff \in\qff G$\sss
on\sss $X\qff \smallsetminus\qff \{\trf a\trf\}$\nnsp.\oss
The group $G_{\dff a}$\sss is\dss canonically\sss isomorphic\sss to $G$\nnsp.\oss
Clearly,\qss 
$C^{\dff 0}\dff(\trf X\dff,\pff b\trf)
\off =\off
C^{\dff 0}\dff(\trf X\dff,\pff \{\dff a\fff,\qff b\trf\}\trf)
\qff \times\qff
G_{\dff a}$
and\dss hence\sss\vspace{3pt}
\[
\quad
H^{\dff 1}\dff(\trf X\dff,\pff b\qff)
\pff =\off
H^{\dff 1}\dff(\trf X\dff,\pff \{\dff a\fff,\qff b\trf\}\qff)\bigl/\fff G_{\dff a}
\qff
\] 

\vspace{-12pt}\vspace{3pt}
Moreover\halfff,\qss
$H^{\dff 1}\dff(\trf X\dff,\pff \{\dff a\fff,\qff b\trf\}\qff)$
can be identified with\sss the product\sss
$H^{\dff 1}\dff(\trf X\dff,\pff b\qff)
\dff \times\dff
G_{\dff a}$\nsp,\oss
but\sss the identification depends on a choice of\dss
a path $p$ connecting $b$ with $a$\nnsp.\oss
The evaluation of\dss cocycles on\sss the path $p$
defines a map\sss
$\varepsilon_{\fff p}\dff \colon\dff
Z^{\dff 1}\dff(\trf X\dff,\pff b\qff)
\qff \ttoo\qff
G\off =\off G_{\dff a}$\nsp,\oss
which,\oss in\sss turn,\oss
leads\sss to a map
$e_{\fff p}\dff \colon\dff
H^{\dff 1}\dff(\trf X\dff,\pff b\qff)
\qff \ttoo\qff
G_{\dff a}$\nsp.\oss
Together\dss with\dss the quotient\dss map
$\mathfrak{q}_{\dff b}\dff \colon\dff
H^{\dff 1}\dff(\trf X\dff,\pff \{\dff a\fff,\qff b\trf\}\qff)
\qff \ttoo\qff
H^{\dff 1}\dff(\trf X\dff,\pff b\qff)$
the map $e_{\fff p}$ leads\sss to a map\vspace{1.5pt}
\[
\quad
f_{\fff p}\dff \colon\dff
H^{\dff 1}\dff(\trf X\dff,\pff \{\dff a\fff,\qff b\trf\}\qff)
\qff \ttoo\qff
H^{\dff 1}\dff(\trf X\dff,\pff b\qff)
\qff \times\qff
G_{\dff a}
\qff
\]

\vspace{-12pt}\vspace{1.5pt}
depending only on $p$\nnsp.\oss
In order\sss to construct\sss an\sss inverse\sss to\sss $f_{\fff p}$ we need\sss
the following\sss lemma.\oss

\mypar{Lemma.}{section}
\emph{Let\dss $z\qff \in\qff Z^{\dff 1}\dff(\trf X\dff,\pff b\qff)$\nnsp.
The cohomology\dss class\dss in\qss
$H^{\dff 1}\dff(\trf X\dff,\pff \{\dff a\fff,\qff b\trf\}\qff)$ 
of\qss  
$\varepsilon_{\fff p}\dff(\dff z\trf)\qff \bullet\qff z$\qss
depends only\sss on\sss $p$\sss and\dss the cohomology\dss class\dss 
$\hclass{z\fff}\qff \in\qff H^{\dff 1}\dff(\trf X\dff,\pff b\qff)$
of\qss the cocycle\sss $z$\nnsp.\oss}

\proof
Since\sss the subset\sss $\{\dff a\fff,\qff b\trf\}$\sss
is\dss discrete,\oss the cocycle $z$\nnsp,\oss
and\dss hence also\sss $\varepsilon_{\fff p}\dff(\dff z\trf)\qff \bullet\qff z$\nnsp,\oss 
automatically\sss belongs\sss to
$Z^{\dff 1}\dff(\trf X\dff,\pff \{\dff a\fff,\qff b\trf\}\qff)$\nnsp.\oss
Suppose\sss that\sss $c\qff \in\qff C^{\dff 0}\dff(\trf X\dff,\pff b\trf)$
and\dss let $w\off =\off c\dff \bullet\dff z$\nnsp.\oss
Then\vspace{3pt}
\[
\quad
\varepsilon_{\fff p}\dff(\dff w\trf)
\off =\off
w\trf(\dff p\trf)
\off =\off
z\dff(\trf p\trf)\dff \cdot\dff
c\trf(\trf a\trf)^{\dff -\dff 1}
\off =\off
\varepsilon_{\fff p}\dff(\dff z\trf)\dff \cdot\dff
c\trf(\trf a\trf)^{\dff -\dff 1}
\pff.
\]

\vspace{-12pt}\vspace{3pt}
By\sss interpreting\sss these equalities as equalities in $G_{\dff a}$
we see\sss that\vspace{4.5pt}
\[
\quad
\varepsilon_{\fff p}\dff(\dff w\trf)\dff \bullet\dff w
\off =\off
\bigl(\qff
\varepsilon_{\fff p}\dff(\dff z\trf)\dff \cdot\trf c\trf(\trf a\trf)^{\dff -\dff 1}
\qff\bigr)
\dff \bullet\dff
\bigl(\qff
c\dff \bullet\dff z
\qff\bigr)
\off =\off
\bigl(\qff
\varepsilon_{\dff p}\dff(\dff z\trf)\dff \cdot\trf 
c\trf(\trf a\trf)^{\dff -\dff 1}
\dff \cdot\dff
c
\qff\bigr)
\dff \bullet\dff
z
\pff.
\]

\vspace{-12pt}\vspace{4.5pt}
Clearly,\oss the $0$\dnsp-cochain
$c\trf(\trf a\trf)^{\dff -\dff 1}\dff \cdot\dff c
\pff \in\qff 
C^{\dff 0}\dff(\trf X\fff,\pff b\trf)$\sss
is\dss equal\sss to $1$ at $a$\nnsp,\oss
and\dss hence
$\varepsilon_{\fff p}\dff(\dff z\trf)\qff \in\qff G_{\dff a}$
and\sss 
$c\trf(\trf a\trf)^{\dff -\dff 1}\dff \cdot\dff c$\sss
commute as elements of\trs the group\dss 
$C^{\dff 0}\dff(\trf X\dff,\pff b\qff)$\nnsp.\oss 
It\sss follows\sss that\vspace{4.5pt}
\[
\quad
\varepsilon_{\fff p}\dff(\dff w\trf)\dff \bullet\dff w
\off =\off
\bigl(\qff 
c\trf(\trf a\trf)^{\dff -\dff 1}
\cdot\dff
c
\dff \cdot\trf
\varepsilon_{\fff p}\dff(\dff z\trf)
\qff\bigr)
\dff \bullet\dff
z
\off =\off
\bigl(\qff
c\trf(\trf a\trf)^{\dff -\dff 1}
\cdot\dff 
c
\qff\bigr)
\dff \bullet\dff
\bigl(\qff
\varepsilon_{\fff p}\dff(\dff z\trf)\dff \bullet\dff z
\qff\bigr)
\pff,
\]

\vspace{-12pt}\vspace{4.5pt}
and\dss hence\sss
$\varepsilon_{\fff p}\dff(\dff w\trf)\dff \bullet\dff w$
and\sss
$\varepsilon_{\dff s}\dff(\dff z\trf)\dff \bullet\dff z$\dss
belong\sss to\sss the same cohomology class.\oss  \eproof

\myuppar{The\sss inverse\sss of $f_{\fff p}$\nsp.}
Lemma\qss \ref{section}\qss implies\sss that\dss the map\dss
$z\off \longmapsto\off \varepsilon_{\fff p}\dff(\dff z\trf)\dff \bullet\dff z$\dss
leads\sss to a map\vspace{3pt}
\[
\quad
\eta_{\trf p}\dff \colon\dff
H^{\dff 1}\dff(\trf X\dff,\pff b\qff)
\qff \ttoo\qff
H^{\dff 1}\dff(\trf X\dff,\pff \{\dff a\fff,\qff b\trf\}\qff)
\qff.
\]

\vspace{-12pt}\vspace{3pt}
Clearly\halfff,\oss 
for every\dss 
$z\qff \in\qff Z^{\dff 1}\dff(\trf X\dff,\pff b\qff)$\dss
the cohomology\sss classes\dss in
$H^{\dff 1}\dff(\trf X\dff,\pff b\qff)$
of\trs $z$\dss and\sss
$\varepsilon_{\dff p}\dff(\dff z\trf)\qff \bullet\qff z$\dss
are equal.\oss
Hence\dss $\eta_{\trf p}$\dss is\dss a section of\trs the map\dss
$\mathfrak{q}_{\dff b}\dff \colon\dff
H^{\dff 1}\dff(\trf X\dff,\pff \{\dff a\fff,\qff b\trf\}\qff)
\qff \ttoo\qff
H^{\dff 1}\dff(\trf X\dff,\pff b\qff)$\nnsp.\oss
An immediate verification shows\sss that\dss
$\varepsilon_{\fff p}\dff(\qff \varepsilon_{\fff p}\dff(\dff z\trf)\dff \bullet\dff z\trf)
\off =\off 
1$\nsp,\oss
and\dss hence\sss $\eta_{\trf p}$ maps\dss
$H^{\dff 1}\dff(\trf X\dff,\pff b\qff)$\dss
bijectively\dss to\sss 
$e_{\fff p}^{\dff -\dff 1}\dff(\dff 1\dff)$\nnsp.\oss
The map $\varepsilon_{\fff p}$\sss is\dss $G_{\dff a}$\dnsp-equivariant\sss
in\sss the sense\sss that\sss
$\varepsilon_{\fff p}\dff(\dff c\dff \bullet\dff u\trf)
\off =\off
\varepsilon_{\fff p}\dff(\dff u\trf)
\dff \cdot\dff
c\trf(\dff a\trf)^{\dff -\dff 1}$\sss
for every\sss $u\qff \in\qff Z^{\dff 1}\dff(\trf X\dff,\pff b\qff)$\nnsp,\dss
$c\qff \in\qff G_{\dff a}$\nsp.\oss
By combining\sss these observations
we see\sss that\sss the map\vspace{3pt}
\[
\quad
g_{\dff p}\dff \colon\dff
H^{\dff 1}\dff(\trf X\dff,\pff b\qff)
\qff \times\qff
G_{\dff a}
\qff \ttoo\qff
H^{\dff 1}\dff(\trf X\dff,\pff \{\dff a\fff,\qff b\trf\}\qff)
\qff,
\]

\vspace{-12pt}\vspace{3pt}
defined\dss by\dss the rule\dss
$g_{\dff p}\dff \colon\dff
(\trf \alpha\dff,\qff c\trf)
\off \longmapsto\off
c^{\dff -\dff 1}\qff \bullet\qff \eta_{\trf p}\dff(\dff \alpha\trf)$\nnsp,\oss
is\dss a\dss bijection.\oss
Another\dss immediate verification shows\sss that\sss $g_{\dff p}$\sss
is\dss the inverse of\dss $f_{\dff p}$\nsp.\oss

\mypar{Theorem.}{svk-olum-relative}
\emph{Suppose\sss that\trs $U\fff,\pff V\qff \subset\qff X$\dss
are\sss two path-connected open sets such\dss that\dss $U\qff \cap\qff V$\dss
consists of\trs two path-connected\sss components\dss
$A$\dss and\dss $B$\nnsp.\oss
If\qss $a\qff \in\qff A$\dss and\dss $b\qff \in\qff B$\nnsp,\oss
then\dss the square}\vspace{3pt}
\[
\quad
\begin{tikzcd}[column sep=large, row sep=normal]
& 
H^{\dff 1}\dff(\trf U\fff,\pff \{\dff a\fff,\qff b\trf\}\qff)
\arrow[rd, near start, "{\dis r_{\trf U\dff \cap\dff V}}"]
&\\
H^{\dff 1}\dff(\trf X\fff,\pff \{\dff a\fff,\qff b\trf\}\qff) 
\arrow[ru, "{\dis r_{\trf U}}"]
\arrow[rd, "{\dis r_{\trf V}}"'] 
&
&
H^{\dff 1}\dff(\trf U\dff \cap\dff V\fff,\pff \{\dff a\fff,\qff b\trf\}\qff)\pff,\\
&
H^{\dff 1}\dff(\trf V\fff,\pff \{\dff a\fff,\qff b\trf\}\qff)
\arrow[ru, near start, "{\dis r_{\trf U\dff \cap\dff V}}"']
&
\end{tikzcd}
\]

\vspace{-9pt}
\emph{where all\dss maps are\sss the restriction maps,\oss
is\dss commutative and cartesian.\oss
There\dss is\dss a\sss canonical\dss bijection\dss between\dss
$H^{\dff 1}\dff(\trf U\dff \cap\dff V\fff,\pff \{\dff a\fff,\qff b\trf\}\qff)$\dss
and\dss
$H^{\dff 1}\dff(\trf A\fff,\pff a\qff)
\dff \times\dff
H^{\dff 1}\dff(\trf B\fff,\pff b\qff)$\nnsp.\oss}

\proof
The second\sss statement\dss is\dss trivial.\oss
The proof\dss of\trs the first\sss one\dss is\dss almost\dss the same as\sss the proof\dss
of\qss Theorem\qss \ref{svk-olum}.\oss
Only\dss the proof\dss of\dss injectivity\dss
used\dss the assumption\dss that\trs $U\qff \cap\qff V$\dss is\dss path-connected.\oss
It\dss was used\dss to ensure\sss that\sss for every\dss 
$x\qff \in\qff U\qff \cap\qff V$\dss there exists\sss
$p\qff \in\qff P\trf(\trf U\qff \cap\qff V\trf)$\dss
connecting\sss
$p\dff(\dff 0\dff)\off =\off b$\sss
with\sss
$p\dff(\dff 1\dff)\off =\off x$\nnsp.\oss
The condition $p\dff(\dff 0\dff)\off =\off b$
was needed\dss to ensure\sss that\sss
$c\trf(\trf p\dff(\dff 0\dff)\trf)\off =\off 1$\sss
for every $c\qff \in\qff C^{\dff 0}\dff(\trf U\fff,\pff b\dff)$\nnsp.\oss
In\dss the present\sss situation\dss we can simply\dss replace\sss this condition\dss
by\dss
$p\dff(\dff 0\dff)\off =\off a$\dss
or\dss
$p\dff(\dff 0\dff)\off =\off b$\nnsp.\oss  \eproof

\mypar{Theorem.}{pseudo-circle}
\emph{Suppose\sss that\trs $X\off =\off U\qff \cup\qff V$\dnsp,\oss
where\dss
$U\fff,\pff V$\dss
are\sss simply\sss connected open sets such\dss that\dss 
$b\qff \in\qff U\qff \cap\qff V$\dss
and\qss $U\qff \cap\qff V$\dss
has\dss two path-connected components.\oss  
Then\dss $\pi_{\dff 1}\dff(\trf X\fff,\pff b\trf)$\dss is\dss a\sss free group
with one generator\halfff,\oss i.e.\qss is\dss isomorphic\sss to\dss $\zzz$\nnsp.\oss}

\proof
Let\dss $B$\sss be\sss the component\sss of\trs the intersection\dss $U\qff \cap\qff V$\dss
containing\sss $b$\nnsp,\oss 
let\dss $A$\sss be\sss the other component\halfff,\oss
and\dss let\trs $a\qff \in\qff A$\nnsp.\oss
Suppose\sss that\dss $W\off =\off U$\sss or\dss $V$\dnsp.\oss
Since\sss $W$\sss is\dss simply-con\-nect\-ed,\oss
$H^{\dff 1}\dff(\trf W\dff,\pff b\qff)$ consists of\dss
one element\halfff,\oss namely\halfff,\oss
the cohomology\sss class of\trs the cocycle 
$\mathbb{1}$ equal\dss to $1$ on every\dss path.\oss
It\dss follows\dss that\sss 
$H^{\dff 1}\dff(\trf W\dff,\pff \{\dff a\fff,\qff b\trf\}\qff)$
can\sss be identified\dss with\dss 
$G_{\dff a}$\nsp.\oss
Under\dss this identification\dss $g\qff \in\qff G_{\dff a}$\dss
corresponds\sss to\sss the cohomology\sss class\sss $h_{\dff g}$\sss of\dss the cocycle\dss 
$g\dff \bullet\dff \mathbb{1}$\nnsp.\oss
In\dss particular\halfff,\oss this identification do not\sss depend on\dss
the choice of\dss $p$\nnsp.\oss
Clearly\halfff,\pss
$g\dff \bullet\dff \mathbb{1}\pff(\trf q\trf)
\off =\off
1$\dss
if\sss $q$\sss is\dss a\sss loop based at\dss $a$\dss or\dss $b$\nnsp.\oss
Therefore\trs Lemma\qss \ref{hom-cohomology}\qss implies\sss that\dss
the images of\trs  
$h_{\dff g}$\dss
in\dss $H^{\dff 1}\dff(\trf A\fff,\pff a\trf)$\dss
and\dss $H^{\dff 1}\dff(\trf B\fff,\pff b\trf)$\dss
are\sss trivial\sss cohomology\sss classes for every\dss $g\qff \in\qff G_{\dff a}$\nnsp.\oss

Now\qss Theorem\qss \ref{svk-olum-relative}\qss implies\dss
that\dss there\dss is\dss a canonical\dss
bijection\dss between\sss
$H^{\dff 1}\dff(\trf X\fff,\pff \{\dff a\fff,\qff b\trf\}\trf)$ 
and\dss $G_{\dff a}\dff \times\dff G_{\dff a}$\nsp.\oss
The group\dss $G_{\dff a}$\dss acts on\sss
$H^{\dff 1}\dff(\trf X\fff,\pff \{\dff a\fff,\qff b\trf\}\trf)$\nnsp.\oss
In\dss terms of\trs $G_{\dff a}\dff \times\dff G_{\dff a}$\dss
this action\sss is\dss the di\-ag\-o\-nal\sss action\dss
$g\dff \bullet\dff (\trf h\fff,\pff k\trf)
\off =\off
(\trf h\dff \cdot\dff g^{\dff -\dff 1}\fff,\pff k\dff \cdot\dff g^{\dff -\dff 1}\trf)$\nnsp.\oss
Therefore\sss the map\dss
$(\trf h\fff,\qff k\trf)
\qff \longmapsto\qff
h\fff \cdot\trf k^{\dff -\dff 1}$\dss is\dss a\sss bijection\qss
$G_{\dff a}\dff \times\dff G_{\dff a}\bigl/\halfff G_{\dff a}
\qff \ttoo\qff
G_{\dff a}$\nsp.\oss
Since\sss $G_{\dff a}$\sss is\dss canonically\dss isomorphic\sss to\sss $G$\sss
and\sss $H^{\dff 1}\dff(\trf X\fff,\pff b\trf)$
is\dss equal\dss to\sss the quotient\sss of\dss
$H^{\dff 1}\dff(\trf X\fff,\pff \{\dff a\fff,\qff b\trf\}\trf)$
by\dss the action of\sss $G_{\dff a}$\nsp,\oss
it\dss follows\dss 
that\dss there\dss is\dss a canonical\dss bijection\dss between
$H^{\dff 1}\dff(\trf X\fff,\pff b\trf)$ and $G$\nnsp.\oss
In\sss view of\qss Section\qss \ref{cochains}\qss
this implies\sss that\dss there\dss is\dss
a\sss canonical\dss bijection\dss between
$\hom (\trf \pi_{\dff 1}\dff(\trf X\fff,\qff b\trf)\fff,\pff G\trf)$
and\dss $G$\nnsp.\oss
By\dss the abstract\dss nonsense,\oss
this means\sss that\dss
$\pi_{\dff 1}\dff(\trf X\fff,\qff b\trf)$\dss
is\dss a\sss free\sss group with one\sss generator\sss
and\dss hence\dss is\dss isomorphic\sss to\sss $\zzz$\nnsp.\oss  \eproof

\myuppar{Changing\dss the base point\halfff.}
Suppose\sss that\dss 
$a\qff \in\qff X$\nnsp,\pss 
$a\off \neq\off b$\nnsp,\oss
and\dss let\sss $p$\dss be a path such\dss that\trs
$p\dff(\dff 0\dff)\off =\off b$\dss
and\dss
$p\dff(\dff 1\dff)\off =\off a$\nnsp.\oss
Recall\dss that\dss the map\dss
$r\off \longmapsto\off p\dff \cdot\dff r\dff \cdot\dff p^{\dff -\dff 1}$\dss
is\dss well\sss defined\dss for\sss loops\sss $r$\sss based at\sss $a$\sss
and after\sss passing\dss to homotopy\sss classes of\dss loops defines
an\dss isomorphism\vspace{3pt}
\[
\quad
\pi\trf(\dff p\trf)\dff \colon\dff
\pi_{\dff 1}\dff(\trf X\dff,\pff a\qff)
\qff \ttoo\qff
\pi_{\dff 1}\dff(\trf X\dff,\pff b\qff)
\qff.
\]

\vspace{-9pt}
Let\dss us\dss consider\dss the composition\dss
$h\trf(\dff p\trf)
\off =\off
\mathfrak{q}_{\dff a}\qff \circ\qff \eta_{\dff p}$,\oss\vspace{4.5pt}
\[
\quad
h\trf(\dff p\trf)
\dff \colon\dff
H^{\dff 1}\dff(\trf X\dff,\pff b\qff)
\qff \ttoo\qff
H^{\dff 1}\dff(\trf X\dff,\pff \{\dff a\fff,\qff b\trf\}\qff)
\qff \ttoo\qff
H^{\dff 1}\dff(\trf X\dff,\pff a\qff)
\qff.
\]

\vspace{-7.5pt}
If\dss $z\qff \in\qff Z^{\dff 1}\dff(\trf X\dff,\pff b\qff)$\dss
and\trs $r$\dss is\dss a\sss loop at\sss $a$\nnsp,\oss
then\dss \vspace{4.5pt}
\[
\quad
\left(\qff\varepsilon_{\dff p}\dff(\dff z\trf)\qff \cdot\qff z \qff\right)\trf(\dff r\trf)
\off =\off
z\trf(\dff p\trf)\dff \cdot\dff 
z\trf(\dff r\trf)\dff \cdot\dff 
z\trf(\dff p\trf)^{\dff -\dff 1}
\off =\off
z\trf\left(\qff  
p\dff \cdot\dff 
r\dff \cdot\dff 
p^{\dff -\dff 1}
\qff\right)
\qff.
\]

\vspace{-7.5pt}
It\dss follows\dss that\dss $h\trf(\dff p\trf)$\dss
is\dss the dual\dss map\sss to\dss $\pi\trf(\dff p\trf)$\dss
in\dss the sense\sss that\dss the diagram\vspace{4.5pt}
\begin{equation*}
\quad
\begin{tikzcd}[column sep=huge, row sep=hugeplus]
H^{\dff 1}\dff(\trf X\fff,\pff b\trf) 
\arrow[r, "{\dis h\trf(\dff p\trf)}"]
\arrow[d] 
&
H^{\dff 1}\dff(\trf X\fff,\pff a\trf)
\arrow[d]\\
\hom (\trf \pi_{\dff 1}\trf(\trf X\fff,\qff b\trf)\fff,\pff G\trf)
\arrow[r, "{\dis \pi\trf(\dff p\trf)^*}"]
&
\hom (\trf \pi_{\dff 1}\trf(\trf X\fff,\qff a\trf)\fff,\pff G\trf)\qff,
\end{tikzcd}
\end{equation*}

\vspace{-7.5pt}
is\dss commutative,\oss
where\dss 
$\pi\trf(\dff p\trf)^*\dff(\dff \varphi\dff)
\off =\off
\varphi\dff \circ\dff \pi\trf(\dff p\trf)$\nnsp.\oss

\myuppar{Changing\dss the base point\sss and\dss the restriction\dss maps.}
Let\sss $a\fff,\qff p$\dss be as above,\oss and\dss 
let\dss
$A\qff \subset\qff X$\dss be a subset\sss such\dss that\dss
$a\qff \in\qff A$\nnsp,\pss
$b\qff \not\in\qff A$\nnsp.\oss
There\dss is\dss a restriction map\vspace{4.5pt}
\[
\quad
\rho_{\dff A}\dff \colon\dff
H^{\dff 1}\dff(\trf X\dff,\pff \{\dff a\fff,\qff b\trf\}\qff)
\qff \ttoo\qff
H^{\dff 1}\dff(\trf A\dff,\pff a\qff)
\qff.
\]

\vspace{-7.5pt}
similar\sss to\sss the restriction maps $r_{\dff U}$\sss from\dss Section\qss \ref{svk-main}.\oss
Clearly\halfff,\pss 
$\rho_{\dff A}
\off =\off
r_{\dff A}\dff \circ\qff \mathfrak{q}_{\dff a}$\sss
and\dss hence\vspace{4.5pt}
\[
\quad
\rho_{\dff A}\dff \circ\qff \eta_{\trf p}
\off =\off
r_{\dff A}\dff \circ\qff \mathfrak{q}_{\dff a}\dff \circ\qff \eta_{\trf p}
\off =\off
r_{\dff A}\dff \circ\qff h\trf(\dff p\trf)
\qff.
\]

\vspace{-7.5pt}
Let\dss
$i\dff \colon\dff
A\qff \ttoo\qff X$\dss
be\sss the inclusion\dss map.\oss
Clearly\halfff,\pss $r_{\dff A}$\dss
is\dss dual\dss to\sss the induced\dss map\dss\vspace{4.5pt}
\[
\quad
i_{\dff *}\dff \colon\dff
\pi_{\dff 1}\trf(\trf A\fff,\qff a\trf)
\qff \ttoo\qff
\pi_{\dff 1}\trf(\trf X\fff,\qff a\trf)
\qff.
\]

\vspace{-7.5pt} 
It\dss follows\dss that\dss
$\rho_{\dff A}\dff \circ\qff \eta_{\trf p}$\dss
is\qss \emph{dual}\qss to\dss
$\pi\trf(\dff p\trf)\dff \circ\dff i_{\dff *}$\dss
in\dss the sense\sss that\dss the diagram\vspace{4.5pt}
\begin{equation*}
\quad
\begin{tikzcd}[column sep=huge, row sep=hugeplus]
H^{\dff 1}\dff(\trf X\fff,\pff b\trf) 
\arrow[r, "{\dis \rho_{\dff A}\dff \circ\qff \eta_{\trf p}}"]
\arrow[d] 
&
H^{\dff 1}\dff(\trf A\fff,\pff a\trf)
\arrow[d]\\
\hom (\trf \pi_{\dff 1}\trf(\trf X\fff,\qff b\trf)\fff,\pff G\trf)
\arrow[r]
&
\hom (\trf \pi_{\dff 1}\trf(\trf A\fff,\qff a\trf)\fff,\pff G\trf)\qff,
\end{tikzcd}
\end{equation*}

\vspace{-7.5pt}
is\dss commutative,\oss
where\sss the lower\dss horizontal\sss arrow\dss is\dss defined\dss by\dss
$\varphi\qff \longmapsto\qff
\varphi\dff \circ\trf \left(\trf \pi\trf(\dff p\trf)\dff \circ\dff i_{\dff *}\dff\right)$\nsp.\oss

\myuppar{Applying\qss Theorem\qss \ref{svk-olum-relative}\qss in\dss the\sss general\sss case.}
Suppose\sss that\dss we are\sss in\dss the situation of\qss
Theorem\qss \ref{svk-olum-relative}.\oss
Let\dss $U\fff'$\dss be a copy\sss of\trs $U$\sss disjoint\dss from\dss $X$\nnsp,\oss
and\dss let\dss $A'\fff,\pff B'\qff \subset\qff U\fff'$\dss be\sss the corresponding\sss
copies of\dss $A\fff,\pff B$\dss respectively\halfff.\oss
Also,\oss let\dss $a'\qff \in\qff A'$\dss be\sss the copy of\dss $a$\nnsp.\oss
One can\dss form a\sss topological\sss space\sss $C$\sss
by\dss identifying\dss $B'$\dss with\dss $B$\dss 
in\dss the union\dss $U\fff'\dff \cup\dff V$\dnsp.\oss
Then\sss the intersection of\sss $U\fff'$ and $V$ in $C$\sss is\dss
equal\dss to $B$ and\dss there\dss is\dss an obvious map\dss 
$\sigma\dff \colon\dff
C\qff \ttoo\qff X$\nnsp.\oss

Our\dss goal\dss is\dss to describe\sss $\pi_{\dff 1}\dff(\trf X\fff,\pff a\trf)$
in\dss terms of\trs $\pi_{\dff 1}\dff(\trf C\fff,\pff a\trf)$\sss
and\sss $\pi_{\dff 1}\dff(\trf A\fff,\pff a\trf)$\nnsp.\oss
Note\sss that\sss $\pi_{\dff 1}\dff(\trf C\fff,\pff a\trf)$\sss
is\dss isomorphic\sss to\sss
$\pi_{\dff 1}\dff(\trf C\fff,\pff b\trf)$\nnsp,\oss
and\dss the\sss latter group\sss 
is\dss determined\dss by\trs
Theorem\qss \ref{svk-universal}\qss as\sss the fundamental\sss group
of\dss the union of\sss $U\fff'$ and $V$\dss
(with identified $B'$ and $B$\nnsp).\oss
We will\sss begin by describing\sss
$H^{\dff 1}\dff(\trf X\dff,\pff \{\trf a\fff,\qff b\trf\}\qff)$\dss
in\sss terms of\dss
$H^{\dff 1}\dff(\trf C\dff,\pff b\trf)$\sss
and\sss
$H^{\dff 1}\dff(\trf A\dff,\pff a\trf)$\nnsp.\oss

Let\dss us\sss choose\sss paths\dss $p\fff,\pff q$\dss connecting\sss $b$\sss
with\sss $a$\sss in\dss $U\fff,\pff V$\dss respectively\halfff.\oss
Let\trs 
$i\dff \colon\dff A\qff \ttoo\qff C$\dss 
be\sss the inclusion\dss map
and\dss let\trs
$\widehat{\imath}\qff \colon\dff A\qff \ttoo\qff C$\dss 
be\sss the composition of\trs the inclusion\dss
$i\fff'\dff \colon\dff 
A'\qff \ttoo\qff C$\dss 
with\dss the\sss tautological\dss homeomorphism\dss
$\iota\dff \colon\dff
A\qff \ttoo\qff A'$\dss
and\dss let\dss $p\fff'$\dss be\sss the copy\sss of\dss $p$\dss in\dss $U\fff'$\nnsp.\oss
The homeomorphism $\iota$ induces a\sss bijection\sss
$\iota^{\fff *}\dff \colon\dff
H^{\dff 1}\dff(\trf A'\dff,\pff a'\qff)
\qff \ttoo\qff
H^{\dff 1}\dff(\trf A\dff,\pff a\qff)$\nnsp.\oss
Let\vspace{4.5pt}
\[
\quad
\rho_{\dff A}\dff \colon\dff
H^{\dff 1}\dff(\trf C\dff,\pff \{\dff a\fff,\qff b\trf\}\qff)
\qff \ttoo\qff
H^{\dff 1}\dff(\trf A\dff,\pff a\qff)
\quad
\mbox{and}\quad
\rho_{\dff A'}\dff \colon\dff
H^{\dff 1}\dff(\trf C\dff,\pff \{\dff a'\fff,\qff b\trf\}\qff)
\qff \ttoo\qff
H^{\dff 1}\dff(\trf A'\dff,\pff a'\qff)
\qff,
\]

\vspace{-12pt}\vspace{7.5pt}
be\sss the obvious restriction maps as above.\oss

\mypar{Theorem.}{two-points-cohomology}
\emph{There\dss is\dss a bijection depending only on\sss $p\fff,\qff q$\sss
between\sss
$H^{\dff 1}\dff(\trf X\dff,\pff \{\dff a\fff,\qff b\trf\}\qff)$\sss
and\sss the set\sss of\dss triples\sss
$(\trf \gamma\fff,\qff h\fff,\qff k\trf)$\sss
with\sss
$\gamma\qff \in\qff H^{\dff 1}\dff(\trf C\fff,\pff b\qff)$
and\sss
$h\fff,\qff k\qff \in\qff G_{\dff a}$\nsp,\oss
such\sss that\sss}\vspace{4.5pt}
\begin{equation}
\label{triple-condition}
\quad
\iota^{\fff *}\dff
\left(\qff
\rho_{\dff A'}\dff \circ\qff \eta_{\dff p\fff'}\trf(\trf \gamma\trf)
\qff\right)
\off =\off
\left(\trf h\dff \cdot\dff k^{\dff -\dff 1}\trf\right)\dff \bullet\trf 
\left(\qff
\rho_{\dff A}\dff \circ\qff \eta_{\dff q}\trf(\qff \gamma\trf)
\qff\right)
\qff.
\end{equation}

\vspace{-7.5pt}
\emph{Moreover\halfff,\oss this bijection\sss turns\sss the action of\dss
$G_{\dff a}$\sss on\sss 
$H^{\dff 1}\dff(\trf X\dff,\pff \{\dff a\fff,\qff b\trf\}\qff)$\sss
into\sss the action}\vspace{4.5pt}
\begin{equation}
\label{triple-action}
\quad
c\dff \bullet\dff (\trf \gamma\fff,\qff h\fff,\qff k\trf)
\off =\off
\left(\trf \gamma\fff,\pff 
h\dff \cdot\dff c\trf(\dff a\trf)^{\dff -\dff 1}\fff,\pff 
k\dff \cdot\dff c\trf(\dff a\trf)^{\dff -\dff 1}\qff\right) 
\qff,
\end{equation}

\vspace{-7.5pt}
\emph{where\sss $c\qff \in\qff G_{\dff a}$\nsp.\oss}

\proof
By\qss Theorem\qss \ref{svk-olum-relative}\qss 
an element\sss of\dss
$H^{\dff 1}\dff(\trf X\dff,\pff \{\dff a\fff,\qff b\trf\}\qff)$\dss
is\dss determined\dss by\dss its\dss images\sss in\vspace{4.5pt}
\begin{equation*}
\quad
H^{\dff 1}\dff(\trf U\fff,\pff \{\dff a\fff,\qff b\trf\}\qff)
\quad\
\mbox{and}\quad\
H^{\dff 1}\dff(\trf V\fff,\pff \{\dff a\fff,\qff b\trf\}\qff)
\qff,
\end{equation*}

\vspace{-7.5pt}
and\sss a pair of\dss cohomology\sss classes in\dss
$H^{\dff 1}\dff(\trf U\fff,\pff \{\dff a\fff,\qff b\trf\}\qff)$
and
$H^{\dff 1}\dff(\trf V\fff,\pff \{\dff a\fff,\qff b\trf\}\qff)$
results\sss from a class in\dss 
$H^{\dff 1}\dff(\trf X\dff,\pff \{\dff a\fff,\qff b\trf\}\qff)$\dss
if\trs and\dss only\trs if\trs their\sss images\sss in\dss
$H^{\dff 1}\dff(\trf U\dff \cap\dff V\fff,\pff \{\dff a\fff,\qff b\trf\}\qff)$\dss
are equal.\oss
Since\dss\vspace{4.5pt}
\[
\quad
H^{\dff 1}\dff(\trf U\dff \cap\dff V\fff,\pff \{\dff a\fff,\qff b\trf\}\qff)
\off =\off
H^{\dff 1}\dff(\trf A\fff,\pff a\qff)
\dff \times\dff
H^{\dff 1}\dff(\trf B\fff,\pff b\qff)
\qff,
\]

\vspace{-7.5pt}
this amounts\sss to\sss
the images in\dss
$H^{\dff 1}\dff(\trf A\fff,\pff a\qff)$\dss
and\dss
$H^{\dff 1}\dff(\trf B\fff,\pff b\qff)$\dss
being\sss equal.\oss
Let\sss us identify\dss the cohomology sets
$H^{\dff 1}\dff(\trf U\fff,\pff \{\dff a\fff,\qff b\trf\}\qff)$\sss
and\sss
$H^{\dff 1}\dff(\trf V\fff,\pff \{\dff a\fff,\qff b\trf\}\qff)$
with\vspace{3pt}
\[
\quad
H^{\dff 1}\dff(\trf U\fff,\pff b\qff)
\qff \times\qff
G_{\dff a}
\quad\
\mbox{and}\quad\
H^{\dff 1}\dff(\trf V\fff,\qff b\qff)
\qff \times\qff
G_{\dff a}
\]

\vspace{-9pt}
by\sss the maps\sss $f_{\fff p}$ and\sss $f_{\fff q}$\sss respectively.\oss
Suppose\sss that\sss\vspace{3pt}
\[
\quad
(\trf \alpha\fff,\qff h\trf)
\qff \in\qff
H^{\dff 1}\dff(\trf U\fff,\pff b\qff)
\qff \times\qff
G_{\dff a} 
\quad\
\mbox{and}\quad\
(\qff \beta\fff,\pff k\trf)
\qff \in\qff
H^{\dff 1}\dff(\trf V\fff,\qff b\qff)
\qff \times\qff
G_{\dff a} 
\qff.
\]

\vspace{-12.25pt}
Then\qss\vspace{-0.25pt}
\[
\quad
g_{\dff p}\dff(\trf \alpha\fff,\qff h\trf)
\off =\off
h^{\dff -\dff 1}\dff \bullet\qff \eta_{\dff p}\trf(\trf \alpha\trf)
\quad\
\mbox{and}\quad\
g_{\dff q}\dff(\trf \beta\fff,\qff k\trf)
\off =\off
k^{\dff -\dff 1}\dff \bullet\qff \eta_{\dff q}\trf(\trf \beta\trf)
\qff.
\]

\vspace{-9pt}
The action of\dss $G_{\dff a}$\dss on cocycles
does not\sss affect\dss their\sss restriction\dss to\sss $B$\nnsp,\oss
and\dss hence\vspace{3pt}
\[
\quad
\rho_{\qff B}\dff\left(\qff 
g_{\dff p}\dff(\trf \alpha\fff,\qff h\trf)
\qff\right)
\off =\off
r_{\qff B}\dff(\trf \alpha\trf)
\quad\
\mbox{and}\quad\
\rho_{\qff B}\dff\left(\qff 
g_{\dff q}\dff(\qff \beta\fff,\qff k\trf)
\qff\right)
\off =\off
r_{\qff B}\dff(\qff \beta\trf)
\qff.
\]

\vspace{-9pt}
On\dss the other\dss hand,\oss
the action of\trs $G_{\dff a}$\dss commutes with\dss 
the restriction\sss to\sss $A$\nnsp,\oss
and\dss hence\vspace{3pt}
\[
\quad
\rho_{\dff A}\dff\left(\qff 
g_{\dff p}\dff(\trf \alpha\fff,\qff h\trf)
\qff\right)
\off =\off
h^{\dff -\dff 1} \bullet\trf 
\rho_{\dff A}\dff\left(\qff 
\eta_{\dff p}\trf(\trf \alpha\trf)
\qff\right)
\quad\
\mbox{and}
\quad\
\rho_{\dff A}\dff\left(\qff 
g_{\dff q}\dff(\qff \beta\fff,\qff k\trf)
\qff\right)
\off =\off
k^{\dff -\dff 1} \bullet\trf 
\rho_{\dff A}\dff\left(\qff 
\eta_{\dff q}\trf(\qff \beta\trf)
\qff\right)
\qff.
\]

\vspace{-9pt}
Hence\sss the images of\sss 
$g_{\dff p}\dff(\trf \alpha\fff,\qff h\trf)$
and\sss
$g_{\dff q}\dff(\trf \beta\fff,\qff k\trf)$
in\sss
$H^{\dff 1}\dff(\trf U\dff \cap\dff V\fff,\pff \{\dff a\fff,\qff b\trf\}\qff)$
are equal\qss if\trs and\dss only\qss if\vspace{3pt}
\[
\quad
r_{\qff B}\dff(\trf \alpha\trf)
\off =\off
r_{\qff B}\dff(\trf \beta\trf)
\quad\
\mbox{and}\quad\
\]

\vspace{-36pt}
\[
\quad
h^{\dff -\dff 1} \bullet\trf 
\left(\qff
\rho_{\dff A}\dff \circ\qff \eta_{\dff p}\trf(\trf \alpha\trf)
\qff\right)
\off =\off
k^{\dff -\dff 1} \bullet\trf 
\left(\qff
\rho_{\dff A}\dff \circ\qff \eta_{\dff q}\trf(\qff \beta\trf)
\qff\right)
\qff.
\]

\vspace{-9pt}
Clearly\halfff,\oss the second condition\dss is\dss equivalent\dss to\vspace{3pt}
\[
\quad
\rho_{\dff A}\dff \circ\qff \eta_{\dff p}\trf(\trf \alpha\trf)
\off =\off
\left(\trf h\dff \cdot\dff k^{\dff -\dff 1}\trf\right)\dff \bullet\trf 
\left(\qff
\rho_{\dff A}\dff \circ\qff \eta_{\dff q}\trf(\qff \beta\trf)
\qff\right)
\qff.
\]

\vspace{-9pt}
By\qss Theorem\qss \ref{svk-olum}\fff\qss
the condition\dss 
$r_{\qff B}\dff(\trf \alpha\trf)
\off =\off
r_{\qff B}\dff(\trf \beta\trf)$\dss
is\dss equivalent\dss to\sss the existence of\dss
a class\dss
$\gamma\qff \in\qff H^{\dff 1}\dff(\trf C\fff,\pff b\qff)$\dss
such\dss that\dss the\sss restriction of\dss $\gamma$\dss to\dss $U\fff'$\dss
is\dss the copy\sss $\alpha'$\sss of\dss $\alpha$\nnsp,\oss
and\dss the\sss restriction of\dss $\gamma$\dss to\dss $V$\trs is\trs $\beta$\nnsp.\oss
Moreover\halfff,\oss when $\gamma$ exists,\oss 
it\dss is\dss unique.\oss
In\dss terms of\dss $\gamma$\dss and\dss the bijection\sss $\iota^{\fff *}$\sss
the second condition\sss takes\sss the form\qss (\ref{triple-condition}).\oss
This proves\sss the first\sss statement\sss of\dss the\sss theorem.\oss

In order\sss to prove\sss the second statement,\oss
note\sss that,\oss similarly\sss to $\varepsilon_{\fff p}$\nsp,\oss 
the map $e_{\fff p}$ is\dss $G_{\dff a}$\dnsp-equivariant\sss
in\sss the sense\sss that\sss
$e_{\fff p}\dff(\dff c\dff \bullet\dff u\trf)
\off =\off
e_{\fff p}\dff(\dff u\trf)
\dff \cdot\dff
c\trf(\dff a\trf)^{\dff -\dff 1}$\sss
for every\sss 
$u\qff \in\qff 
H^{\dff 1}\dff(\trf U\fff,\pff \{\dff a\fff,\qff b\trf\}\qff)$\nnsp,\dss
$c\qff \in\qff G_{\dff a}$\nsp.\oss
At\sss the same\sss time,\oss clearly,\dss
$\mathfrak{q}_{\dff b}\dff(\dff c\dff \bullet\dff u\trf)
\off =\off
\mathfrak{q}_{\dff b}\dff(\dff u\trf)$\sss
for every\sss 
$u\qff \in\qff 
H^{\dff 1}\dff(\trf U\fff,\pff \{\dff a\fff,\qff b\trf\}\qff)$\sss
and\sss $c\qff \in\qff G_{\dff a}$\nsp.\oss
It\sss follows\sss that\sss $f_{\fff p}$\sss turns\sss the action of\dss $G_{\dff a}$\sss
on\sss $H^{\dff 1}\dff(\trf U\fff,\pff \{\dff a\fff,\qff b\trf\}\qff)$
into\sss the action\sss\vspace{3pt}
\[
\quad
c\dff \bullet\dff (\trf \alpha\fff,\qff h\trf)
\qff \longmapsto\qff
\left(\trf \alpha\fff,\qff h\dff \cdot\dff c\trf(\dff a\trf)^{\dff -\dff 1}\trf\right)
\]

\vspace{-9pt}
The map $f_{\fff q}$ has a similar\sss property.\oss
The second statement\sss follows.\oss  \eproof

\mypar{Corollary.}{two-point-quot}
\emph{There\dss is\dss a bijection depending only on\sss $p\fff,\qff q$\sss
between\sss
$H^{\dff 1}\dff(\trf X\dff,\pff b\qff)$\sss
and\sss the set\sss of\dss pairs\sss
$(\trf \gamma\fff,\qff g\trf)$\sss
with\sss
$\gamma\qff \in\qff H^{\dff 1}\dff(\trf C\fff,\pff b\qff)$
and\sss
$g\qff \in\qff G_{\dff a}$\nsp,\oss
such\sss that\sss}\vspace{3pt}
\begin{equation}
\label{pair-condition}
\quad
\iota^{\fff *}\dff
\left(\qff
\rho_{\dff A'}\dff \circ\qff \eta_{\dff p\fff'}\trf(\trf \gamma\trf)
\qff\right)
\off =\off
g\dff \bullet\trf 
\left(\qff
\rho_{\dff A}\dff \circ\qff \eta_{\dff q}\trf(\qff \gamma\trf)
\qff\right)
\qff.
\end{equation}

\vspace{-9pt}
\proof
It\dss is\dss sufficient\sss to notice\sss that\sss 
the action\qss (\ref{triple-action})\qss leaves\sss the product\dss
$h\dff \cdot\dff k^{\dff -\dff 1}$\dss invariant.\oss  \eproof

\myuppar{The fundamental\dss groups.}
Let $\tau\qff \in\qff \pi_{\dff 1}\dff(\trf X\fff,\qff b\trf)$ be\sss
the homotopy\sss class of\dss the\sss loop $p\dff \cdot\dff q^{\dff -\dff 1}$\dnsp.\oss
Let\vspace{3pt}
\[
\quad
\theta
\off =\off
\pi\trf(\dff q\trf)\trf \circ\qff i_{\dff *}
\dff \colon\dff
\pi_{\dff 1}\dff(\trf A\fff,\qff a\trf)
\off \ttoo\off
\pi_{\dff 1}\dff(\trf C\fff,\qff b\trf)
\quad
\mbox{and}\quad
\]

\vspace{-36pt}
\[
\quad
\widehat{\theta}
\off =\off
\pi\trf(\dff p'\trf)\trf \circ\qff \widehat{\imath}_{\dff *}
\dff \colon\dff
\pi_{\dff 1}\dff(\trf A\fff,\qff a\trf)
\off \ttoo\off
\pi_{\dff 1}\dff(\trf C\fff,\qff b\trf)
\qff.
\]

\vspace{-9pt}
\mypar{Theorem.}{open-van-kampen}
\emph{Let\qss $F$ be\sss the free\sss group\sss with\dss
one\sss generator\dss $t$\nnsp.\oss
The\sss fundamental\dss group\dss
$\pi_{\dff 1}\dff(\trf X\fff,\qff b\trf)$\dss
is\dss isomorphic\sss to\sss the\dss free product\trs
$\pi_{\dff 1}\dff(\trf C\fff,\qff b\trf)\trf *\qff F$\dss
with\dss the relation}\vspace{3pt}
\[
\quad
\widehat{\theta}\trf(\trf \alpha\trf)
\off =\off
t\pff \theta\trf(\trf \alpha\trf)\qff t^{\dff -\dff 1}
\]

\vspace{-9pt}
\emph{imposed\trs for every\dss 
$\alpha\qff \in\qff \pi_{\dff 1}\dff(\trf A\fff,\qff a\trf)$\nnsp.\oss
The corresponding\dss homomorphism}\vspace{3pt}
\[
\quad
\pi_{\dff 1}\dff(\trf C\fff,\qff b\trf)\trf *\qff F
\off \ttoo\off
\pi_{\dff 1}\dff(\trf X\fff,\qff b\trf)
\]

\vspace{-9pt}
\emph{is\dss equal\dss to\dss
$\sigma_{*}\dff \colon\dff
\pi_{\dff 1}\dff(\trf C\fff,\qff b\trf)
\qff \ttoo\qff
\pi_{\dff 1}\dff(\trf X\fff,\qff b\trf)$\dss
on\qss
$\pi_{\dff 1}\dff(\trf C\fff,\qff b\trf)$\dss
and\dss maps\dss $t$\dss to\dss $\tau$\nnsp.\oss}

\proof
Let\sss use\dss Theorem\qss \ref{independence-of-paths}\qss
to pass\sss from\sss $H^{\dff 1}\dff(\trf X\dff,\pff b\qff)$\sss to\sss
$\hom (\trf \pi_{\dff 1}\dff(\trf X\fff,\qff b\trf)\fff,\pff G\trf)$\nnsp.\oss
We can\sss identify\sss this set\sss of\dss homomorphisms 
with\sss the set\sss of\dss pairs\sss
$(\trf \varphi\fff,\qff g\trf)$\sss such\sss that\sss $\varphi$\sss
is\dss a homomorphism\sss
$\pi_{\dff 1}\dff(\trf C\fff,\qff b\trf)\qff \ttoo\qff G$\nnsp,\dss
$g\qff \in\qff G_{\dff a}\off =\off G$\nnsp,\oss
and\sss the cohomology class $\gamma$ corresponding\sss to $\varphi$\sss
together with $g$ satisfies\qss (\ref{pair-condition}).\oss
Let\sss us reformulate\qss (\ref{pair-condition})\qss 
in\sss terms of\trs homomorphisms $\varphi$\nnsp.\oss
As we saw,\qss
$\rho_{\dff A}\dff \circ\qff \eta_{\trf q}$\dss
is\dss dual\dss to\dss
$\theta
\off =\off
\pi\trf(\dff q\trf)\dff \circ\dff i_{\dff *}$\sss
and\sss
$\rho_{\dff A'}\dff \circ\qff \eta_{\trf p'}$\dss
is\dss dual\dss to\dss
$\widehat{\theta}
\off =\off
\pi\trf(\dff p'\trf)\trf \circ\qff \widehat{\imath}_{\dff *}$\nsp.\oss
The action of\dss the group\sss $G\off =\off G_{\dff a}$ on\sss
$H^{\dff 1}\dff(\trf A\dff,\pff a\qff)$\sss
corresponds\sss to\sss the action of\sss $G$ on\sss
$\hom (\trf \pi_{\dff 1}\dff(\trf A\fff,\qff a\trf)\fff,\pff G\trf)$\sss
by conjugation,\oss namely\sss to\sss the action\sss
$(\trf g\fff,\qff \psi\trf)
\qff \longmapsto\qff
g\dff \psi\dff g^{\dff -\dff 1}$\dnsp,\oss
where $g\dff \psi\dff g^{\dff -\dff 1}$\sss 
is\dss the homomorphism\sss
$r
\qff \longmapsto\qff 
g\dff \cdot\dff \psi\dff(\dff r\trf)\dff \cdot\dff g^{\dff -\dff 1}$\dnsp.\oss
It\sss follows\sss that\qss (\ref{pair-condition})\qss holds for\sss
$(\trf \varphi\fff,\qff g\trf)$\sss
if\dss and\dss only\trs if\vspace{3pt}
\begin{equation}
\label{hom-condition-orbits}
\quad
\varphi\dff \circ\qff \widehat{\theta}
\off =\off
g\dff 
\left(\qff
\varphi\dff \circ\dff \theta
\qff\right)\dff g^{\dff -\dff 1}
\qff.
\end{equation}

\vspace{-9pt}
In\sss turn,\oss this implies\sss that\dss the group\dss
$\pi_{\dff 1}\dff(\trf X\fff,\qff b\trf)$\dss
has\sss the same universal\dss property\sss as\sss the quotient\dss
group of\dss 
$\pi_{\dff 1}\dff(\trf C\fff,\qff b\trf)\trf *\qff F$\dss
described\dss in\dss the\sss theorem.\oss
Hence\dss
$\pi_{\dff 1}\dff(\trf X\fff,\qff b\trf)$\dss
is\dss isomorphic\sss to\sss this quotient\halfff.\oss
It\dss remains\sss to check\sss the claim about\dss the image of\dss $t$\dss
in\dss $\pi_{\dff 1}\dff(\trf X\fff,\qff b\trf)$\nnsp.\oss
Since\sss the maps\dss $g_{\dff p}\dff,\pff g_{\dff q}$\dss
are\sss the inverses of\trs the maps\dss $f_{\dff p}\dff,\pff f_{\dff q}$\dss
respectively\halfff,\oss
we see\sss that\vspace{3pt}
\[
\quad
e_{\dff p}\trf\left(\qff
g_{\dff p}\dff(\trf \alpha\fff,\qff h\trf)
\qff\right)
\off =\off
h
\quad\
\mbox{and}\quad\
e_{\dff q}\trf\left(\qff
g_{\dff q}\dff(\trf \beta\fff,\qff k\trf)
\qff\right)
\off =\off
k
\qff.
\]

\vspace{-9pt}
It\dss follows\dss that\dss the cohomology\dss class\sss in\dss
$H^{\dff 1}\dff(\trf X\dff,\pff \{\dff a\fff,\qff b\trf\}\qff)$\dss
defined\dss by\dss $g_{\dff p}\dff(\trf \alpha\fff,\qff h\trf)$\dss
and\dss $g_{\dff q}\dff(\trf \beta\fff,\qff k\trf)$
takes\sss the value\dss $h\dff \cdot\dff k^{\dff -\dff 1}$\dss
on\dss the path\dss $p\dff \cdot\dff q^{\dff -\dff 1}$\dnsp.\oss
The same\dss is\dss true for\dss the image of\trs this cohomology\sss
class\sss in $H^{\dff 1}\dff(\trf X\dff,\pff b\qff)$\nnsp.\oss
In\dss terms of\trs the corresponding\dss homomorphism\dss
$\psi\dff \colon\dff
\pi_{\dff 1}\dff(\trf X\fff,\qff b\trf)
\qff \ttoo\qff
G$\dss
this means\sss that\dss
$\psi\dff(\dff \tau\trf)\off =\off h\dff \cdot\dff k^{\dff -\dff 1}$\dnsp.\oss
Now\dss the claim about\dss the image of\dss $t$\dss follows\sss from\dss
the abstract\dss nonsense.\oss
One can also apply\dss the universal\dss property\sss of\dss
$\pi_{\dff 1}\dff(\trf X\fff,\qff b\trf)$\dss to\dss
$G\off =\off \pi_{\dff 1}\dff(\trf X\fff,\qff b\trf)$\dss and\dss the identity\dss
homomorphism.\oss
We\sss leave\sss the details\sss to\sss the reader\halfff.\oss  \eproof

\myuppar{The case of\dss general\sss discrete subsets $Y$\dnsp.}
Now we are going\sss to discuss\sss
the situation when\sss the intersection $U\dff \cap\dff V$\sss
consists of\dss several\sss path-connected components.\oss
This requires some preliminary discussion of\dss
$H^{\dff 1}\dff(\trf X\dff,\pff Y\trf)$
with discrete subset\sss $Y$ similar\sss to\sss the discussion
at\sss the beginning of\dss this section.\oss
So,\oss let\sss us assume\sss that $Y$ is\dss discrete\qss
(actually,\oss it\dss is\dss sufficient\sss to assume\sss
that\sss every map\sss $[\dff 0\fff,\qff 1\dff]\qff \ttoo\qff Y$\sss is\dss constant).\oss
Let\dss us\dss fix\dss for every\dss $y\qff \in\qff Y$\sss
a\dss path\sss $s_{\dff y}$\sss such\dss that\trs
$s_{\dff y}\trf(\dff 0\dff)\off =\off b$\dss
and\dss
$s_{\dff y}\trf(\dff 1\dff)\off =\off y$\nnsp.\oss
We will\sss assume\sss that\dss $s_{\dff b}$\sss is\dss
constant\halfff.\oss
The eval\-u\-a\-tion of\dss cocycles on\dss paths\sss $s_{\dff y}$\sss
defines a map\sss\vspace{4.5pt}
\[
\quad
\varepsilon_{\dff s}\dff \colon\dff
Z^{\dff 1}\dff(\trf X\dff,\pff b\qff)
\qff \ttoo\qff
C^{\dff 0}\dff(\trf Y\fff,\pff b\trf)
\qff,
\]

\vspace{-7.5pt}
which,\oss in\dss turn,\oss leads\sss another evaluation\dss map\sss 
$e_{\fff s}\dff \colon\dff
H^{\dff 1}\dff(\trf X\dff,\pff Y\qff)
\qff \ttoo\qff
C^{\dff 0}\dff(\trf Y\fff,\pff b\trf)$\nnsp.\oss
Together\dss with\dss the quotient\dss map\dss 
$\mathfrak{q}_{\dff b}\dff \colon\dff
H^{\dff 1}\dff(\trf X\dff,\pff Y\qff)
\qff \ttoo\qff
H^{\dff 1}\dff(\trf X\dff,\pff b\qff)$\dss
the map\dss $e_{\fff s}$\dss leads\sss to a map\vspace{4.5pt}
\[
\quad
f_{\fff s}\dff \colon\dff
H^{\dff 1}\dff(\trf X\dff,\pff Y\qff)
\qff \ttoo\qff
H^{\dff 1}\dff(\trf X\dff,\pff b\qff)
\qff \times\qff
C^{\dff 0}\dff(\trf Y\fff,\pff b\trf)
\qff.
\]

\vspace{-7.5pt}
In order\sss to prove\sss that\sss $f_{\fff s}$\sss is\dss a bijection,\oss
we need\sss the following\sss lemma,\oss
in which we identify\sss $C^{\dff 0}\dff(\trf Y\fff,\pff b\trf)$
with\sss the subgroup of\dss $C^{\dff 0}\dff(\trf X\fff,\pff b\trf)$
consisting of\sss $0$\dnsp-cochains equal\dss to $1$ on $X\dff \smallsetminus\dff Y$\dnsp.\oss

\mypar{Lemma.}{section-general}
\emph{Let\dss $z\qff \in\qff Z^{\dff 1}\dff(\trf X\dff,\pff b\qff)$\nnsp.
The cohomology\dss class\dss in\qss
$H^{\dff 1}\dff(\trf X\dff,\pff Y\qff)$ 
of\qss  
$\varepsilon_{\dff s}\dff(\dff z\trf)\qff \bullet\qff z$\qss
depends only\sss on\sss $s$\sss and\dss the cohomology\dss class\dss 
in\dss
$H^{\dff 1}\dff(\trf X\dff,\pff b\qff)$
of\qss the cocycle\dss $z$\nnsp.\oss}

\proof
The proof\trs is\dss completely similar\sss to\sss 
the proof\dss of\trs Lemma\qss \ref{section}.\oss  \eproof

\myuppar{The\sss inverse\sss of $f_{\fff s}$\nsp.}
Lemma\qss \ref{section-general}\qss implies\sss that\dss the map\dss
$z\off \longmapsto\off \varepsilon_{\fff s}\dff(\dff z\trf)\dff \bullet\dff z$\dss
leads\sss to a map\vspace{4.5pt}
\[
\quad
\eta_{\trf p}\dff \colon\dff
H^{\dff 1}\dff(\trf X\dff,\pff b\qff)
\qff \ttoo\qff
H^{\dff 1}\dff(\trf X\dff,\pff Y\qff)
\qff.
\]

\vspace{-12pt}\vspace{7.5pt}
Similarly\sss to\sss the case of\dss two-points subsets $Y$\dnsp,\oss
the map\vspace{4.5pt}
\[
\quad
g_{\dff s}\dff \colon\dff
H^{\dff 1}\dff(\trf X\dff,\pff b\qff)
\qff \times\qff
C^{\dff 0}\dff(\trf Y\fff,\pff b\trf)
\qff \ttoo\qff
H^{\dff 1}\dff(\trf X\dff,\pff Y\qff)
\qff,
\]

\vspace{-12pt}\vspace{7.5pt}
defined\dss by\dss the rule\dss
$g_{\dff s}\dff \colon\dff
(\trf \alpha\dff,\qff c\trf)
\off \longmapsto\off
c^{\dff -\dff 1}\qff \bullet\qff \eta_{\trf s}\dff(\dff \alpha\trf)$\nnsp,\oss
is\dss a\dss bijection\sss
and\dss is\dss the inverse of\dss $f_{\dff s}$\nsp.\oss
Now we are almost\sss ready\sss to a generalization of\qss 
Theorem\qss \ref{pseudo-circle}.\oss

\mypar{Theorem.}{svk-olum-relative-general}
\emph{Suppose\sss that\trs $U\fff,\pff V\qff \subset\qff X$\dss
are\sss two path-connected open sets such\dss that\dss $U\qff \cap\qff V$\dss
has\dss $m\qff +\qff 1$\dss path-connected components\dss
$A_{\dff 0}\dff,\pff A_{\dff 1}\dff,\pff
\ldots\dff,\pff A_{\dff m}$\nsp,\pss $m\qff \geq\qff 1$\nnsp.\oss
Suppose\sss that\trs $a_{\dff i}\qff \in\qff A_{\dff i}$\dss
for every\dss $i$\dss and\dss let\pss
$Y
\off =\off 
\{\trf a_{\dff 0}\dff,\off a_{\dff 1}\dff,\off \ldots\dff,\off a_{\dff m}\trf\}$\nnsp.\oss
Then\dss the square}\vspace{7pt}
\[
\quad
\begin{tikzcd}[column sep=large, row sep=normal]
& 
H^{\dff 1}\dff(\trf U\fff,\pff Y\qff)
\arrow[rd, near start, "{\dis r_{\trf U\dff \cap\dff V}}"]
&\\
H^{\dff 1}\dff(\trf X\fff,\pff Y\qff) 
\arrow[ru, "{\dis r_{\trf U}}"]
\arrow[rd, "{\dis r_{\trf V}}"'] 
&
&
H^{\dff 1}\dff(\trf U\dff \cap\dff V\fff,\pff Y\qff)\pff,\\
&
H^{\dff 1}\dff(\trf V\fff,\pff Y\qff)
\arrow[ru, near start, "{\dis r_{\trf U\dff \cap\dff V}}"']
&
\end{tikzcd}
\]

\vspace{-5pt}
\emph{where all\dss maps are\sss the restriction maps,\oss
is\dss commutative and cartesian.\oss
There\dss is\dss a\sss canonical\dss bijection\dss between\dss
$H^{\dff 1}\dff(\trf U\dff \cap\dff V\fff,\pff Y\qff)$\dss
and\dss the product\pss
$\prod_{\qff i}\dff H^{\dff 1}\dff(\trf A_{\dff i}\fff,\pff a_{\dff i}\qff)$\nnsp.\oss}

\proof
The proof\trs is\dss a\sss direct\dss generalization of\trs the proof\dss
of\pss Theorem\qss \ref{svk-olum-relative}.\oss  \eproof

\mypar{Theorem.}{pseudo-wedge-of-circles}
\emph{Suppose\sss that\trs $X\off =\off U\qff \cup\qff V$\dnsp,\oss
where\dss
$U\fff,\pff V$\dss
are\sss simply\sss connected open sets such\dss that\dss 
$b\qff \in\qff U\qff \cap\qff V$\dss
and\qss $U\qff \cap\qff V$\dss
has\dss $m\qff +\qff 1$\dss path-connected components,\pss $m\qff \geq\qff 1$\nnsp.\oss 
Then\dss $\pi_{\dff 1}\dff(\trf X\fff,\pff b\trf)$\dss is\dss a\sss free group
with\sss $m$\sss generators.\oss}

\proof
The proof\trs is\dss a\sss direct\dss generalization of\trs the proof\dss
of\pss Theorem\qss \ref{pseudo-circle}.\oss  \eproof

\myuppar{Intersections\sss with\dss more\sss than\dss two components\sss in\sss general.}
Suppose\sss that\dss we are\sss in\dss the situation of\qss
Theorem\qss \ref{svk-olum-relative-general}.\oss
Theorem\qss \ref{open-van-kampen}\qss can be easily\sss extended\dss
to\sss this situation.\oss
In order\dss to stress\sss the analogy\dss with\trs 
Theorem\qss \ref{open-van-kampen},\oss
let\dss us\dss set\dss $B\off =\off A_{\dff 0}$\dss and\dss 
$b\off =\off a_{\dff 0}$\nsp.\oss

To begin\dss with,\oss
let\dss $U\fff'$\dss be a copy\sss of\trs $U$\sss disjoint\dss from\dss $X$\nnsp,\oss
and\dss let\dss $B'\qff \subset\qff U\fff'$\dss be\sss the 
copy\sss of\dss $B$\nnsp.\oss
Let\dss us\dss form a\sss topological\sss space\sss $C$\sss
by\dss identifying\dss $B'$\dss with\dss $B$\dss 
in\dss the union\dss $U\fff'\qff \cup\qff V$\dnsp.\oss
There\dss is\dss an obvious map\dss 
$\sigma\dff \colon\dff
C\qff \ttoo\qff X$\nnsp.\oss
Let\dss $A'_{\dff k}\off \subset\off U\fff'$\qss be\sss the copy\sss
of\dss $A_{\dff k}$\nsp.\oss

Let\sss us\dss choose for each\dss
$k\off =\off 1\fff,\pff 2\fff,\pff \ldots\fff,\pff m$\dss 
paths\dss $p_{\dff k}\dff,\off q_{\dff k}$\dss
connecting\sss $b$\sss with\sss $a_{\dff k}$\sss in\dss $U$\dss and\dss $V$\dss
respectively\halfff.\oss
Then\dss $p_{\dff k}\dff \cdot\dff q_{\dff k}^{\dff -\dff 1}$
are\sss loops in\dss $X$\dss 
based\sss at\dss $b$\nnsp.\oss
Let\dss $\tau_{\dff k}\qff \in\qff \pi_{\dff 1}\dff(\trf X\fff,\qff b\trf)$\dss be\sss
the\sss homotopy\sss class of\trs $p_{\dff k}\dff \cdot\dff q_{\dff k}^{\dff -\dff 1}$\dnsp.\oss
Let\trs 
$i_{\dff k}\dff \colon\dff A_{\dff k}\qff \ttoo\qff C$\dss 
be\sss the inclusion\dss map
and\trs let\dss\vspace{4.5pt}
\[
\quad
\theta_{\dff k}
\off =\off\dff
\pi\trf\left(\qff q_{\dff k} \qff\right)\trf \circ\qff 
\bigl(\qff i_{\dff k} \qff\bigr)_{\dff *}
\qff \colon\qff
\pi_{\dff 1}\qff\bigl(\trf A_{\dff k}\fff,\qff a_{\dff k}\trf\bigr)
\off \ttoo\off
\pi_{\dff 1}\dff(\trf C\fff,\qff b\trf)
\qff.
\]

\vspace{-7.5pt}
Let\trs
$\widehat{\imath}_{\dff k}\qff \colon\dff A_{\dff k}\qff \ttoo\qff C$\dss 
be\sss the composition of\trs the inclusion\dss
$A'_{\dff k}\qff \ttoo\qff C$\dss 
with\dss the\sss tautological\dss homeomorphism\dss
$\iota_{\dff k}\dff \colon\dff
A_{\dff k}\qff \ttoo\qff A'_{\dff k}$\dss
and\dss let\dss $p\fff'_{\dff k}$\dss be\sss 
the copy\sss of\dss $p_{\dff k}$\dss in\dss $U\fff'$\nnsp.\oss
Let\dss\vspace{4.5pt}
\[
\quad
\widehat{\theta}_{\dff k}
\off =\off\dff
\pi\trf\left(\qff p'_{\dff k} \qff\right)
\trf \circ\qff 
\bigl(\qff \widehat{\imath}_{\dff k} \qff\bigr)_{\dff *}
\qff \colon\qff
\pi_{\dff 1}\qff\bigl(\trf A_{\dff k}\fff,\qff a_{\dff k}\trf\bigr)
\off \ttoo\off
\pi_{\dff 1}\dff(\trf C\fff,\qff b\trf)
\qff.
\]

\vspace{-7.5pt}
\mypar{Theorem.}{open-van-kampen-general}
\emph{Let\qss $F$ 
be\sss the free\sss group\sss with\dss
$m$\sss generators\dss $t_{\dff 1}\dff,\pff \ldots\dff,\pff t_{\dff m}$\nsp.\oss
The\sss fundamental\dss group\dss
$\pi_{\dff 1}\dff(\trf X\fff,\qff b\trf)$\dss
is\dss isomorphic\sss to\sss the\dss free product\trs
$\pi_{\dff 1}\dff(\trf C\fff,\qff b\trf)
\trf *\qff 
F$\dss
with\dss the relation}\vspace{4.5pt}
\[
\quad
\widehat{\theta}_{\dff k}\trf(\trf \alpha\trf)
\off =\off
t_{\dff k}\pff \theta_{\dff k}\trf(\trf \alpha\trf)\qff t_{\dff k}^{\dff -\dff 1}
\]

\vspace{-7.5pt}
\emph{imposed\trs for every\dss 
$\alpha
\qff \in\qff 
\pi_{\dff 1}\qff\bigl(\trf A_{\dff k}\fff,\qff a_{\dff k}\trf\bigr)$\dss
and every\dss
$k$\nnsp.\oss
The corresponding\dss homomorphism}\vspace{4.5pt}
\[
\quad
\pi_{\dff 1}\dff(\trf C\fff,\qff b\trf)\trf *\qff F
\off \ttoo\off
\pi_{\dff 1}\dff(\trf X\fff,\qff b\trf)
\]

\vspace{-7.5pt}
\emph{is\dss equal\dss to\dss
$\sigma_{*}\dff \colon\dff
\pi_{\dff 1}\dff(\trf C\fff,\qff b\trf)
\qff \ttoo\qff
\pi_{\dff 1}\dff(\trf X\fff,\qff b\trf)$\dss
on\qss
$\pi_{\dff 1}\dff(\trf C\fff,\qff b\trf)$\dss
and\dss maps\dss $t_{\dff k}$\dss to\dss $\tau_{\dff k}$\nsp.\oss}

\proof
The proof\dss is\dss similar\dss to\sss the proof\dss
of\pss Theorem\qss \ref{open-van-kampen}.\oss 
The main difference\dss is\dss 
the need\dss to use more cumbersome notations.\oss \eproof

\vspace*{\medskipamount}

\begin{flushright}

September\qss 23\fff,\oss 2023
 
https\halfff:/\!/nikolaivivanov.com

E-mail\halfff:\oss nikolai.v.ivanov{\fff}@{\dff}icloud.com,\oss ivanov{\fff}@{\dff}msu.edu

Department\sss of\qss Mathematics,\oss Michigan\sss State\sss University

\end{flushright}

}

\end{document}